\documentclass[11pt]{amsart}

\usepackage[top=2.6cm,
bottom=2.6cm,
left=2.3cm,
right=2.3cm]{geometry}

\usepackage[T1]{fontenc}
\usepackage{amsfonts}
\usepackage{graphicx}
\usepackage{epstopdf}
\usepackage{enumitem}

\newtheorem{theorem}{Theorem}[section]

\newtheorem{proposition}[theorem]{Proposition}

\newtheorem{conjecture}[theorem]{Conjecture}
\newtheorem{corollary}[theorem]{Corollary}

\theoremstyle{definition}

\theoremstyle{remark}
\newtheorem{example}[theorem]{Example}
\newtheorem{remark}[theorem]{Remark}

\usepackage{mathrsfs}

\newcommand{\Yang}{\mathcal{Y}}
\newcommand{\MO}{\mathtt{MO}}

\newcommand{\Z}{{\ensuremath{{\mathbb Z}}}}

\newcommand{\bA}{{\mathbb A}}

\newcommand{\bC}{{\mathbb C}}
\newcommand{\bD}{{\mathbb D}}
\newcommand{\bE}{{\mathbb E}}
\newcommand{\bF}{{\mathbb F}}

\newcommand{\bL}{{\mathbb L}}

\newcommand{\bN}{{\mathbb N}}

\newcommand{\bQ}{{\mathbb Q}}
\newcommand{\bR}{{\mathbb R}}

\newcommand{\bZ}{{\mathbb Z}}
     
\newcommand{\cA}{{\mathcal A}}

\newcommand{\cF}{{\mathcal F}}
\newcommand{\cG}{{\mathcal G}}
\newcommand{\cH}{{\mathcal H}}

\newcommand{\cM}{{\mathcal M}}
\newcommand{\cN}{{\mathcal N}}
\newcommand{\cO}{{\mathcal O}}

\newcommand{\cZ}{{\mathcal Z}}

\newcommand{\fgl}{{\mathfrak {gl}}}

\newcommand{\fh}{{\mathfrak h}}

\newcommand{\fg}{{\mathfrak g}}
\newcommand{\fn}{{\mathfrak n}}

\newcommand{\fw}{{\mathfrak w}}

\newcommand{\fM}{\mathfrak{M}}
\newcommand{\fP}{\mathfrak{P}}
\newcommand{\fS}{\mathfrak{S}}

\newcommand{\sA}{\mathscr{A}}
\newcommand{\sC}{\mathscr{C}}

\DeclareMathOperator{\mW}{W}
\DeclareMathOperator{\simp}{simp}
\DeclareMathOperator{\Coh}{Coh}

\newcommand{\cref}[1]{\ref{#1}}

\makeatletter
\newcommand{\colim@}[2]{%
  \vtop{\m@th\ialign{##\cr
    \hfil$#1\operator@font colim$\hfil\cr
    \noalign{\nointerlineskip\kern1.5\ex@}#2\cr
    \noalign{\nointerlineskip\kern-\ex@}\cr}}%
}
\newcommand{\colim}{%
  \mathop{\mathpalette\colim@{\rightarrowfill@\scriptscriptstyle}}\nmlimits@
}
\renewcommand{\varprojlim}{%
  \mathop{\mathpalette\varlim@{\leftarrowfill@\scriptscriptstyle}}\nmlimits@
}
\renewcommand{\varinjlim}{%
  \mathop{\mathpalette\varlim@{\rightarrowfill@\scriptscriptstyle}}\nmlimits@
}
\makeatother

\newcommand{\1}{\mathbf{1}}

\newcommand{\id}{\mathrm{id}}

\DeclareMathOperator{\Exp}{\mathbf{Exp}}

\DeclareMathOperator{\KK}{K}
\DeclareMathOperator{\eff}{eff}
\DeclareMathOperator{\UEA}{U}

\DeclareMathOperator{\mult}{\mathbf{m}}

\DeclareMathOperator{\Hom}{Hom}
\DeclareMathOperator{\HO}{H}

\DeclareMathOperator{\lmod}{-mod}

\DeclareMathOperator{\IC}{\mathcal{IC}}

\DeclareMathOperator{\Span}{Span}
\DeclareMathOperator{\kac}{\mathtt{a}}

\DeclareMathOperator{\Sym}{\mathrm{Sym}}
\DeclareMathOperator{\GL}{\text{GL}}

\DeclareMathOperator{\opp}{op}

\DeclareMathOperator{\Spec}{Spec}
\DeclareMathOperator{\sst}{-sst}
\DeclareMathOperator{\twisted}{tw}

\DeclareMathOperator{\Gr}{Gr}

\DeclareMathOperator{\frex}{-fr}

\DeclareMathOperator{\pt}{pt}

\DeclareMathOperator{\Ext}{Ext}
\DeclareMathOperator{\Jac}{Jac}

\DeclareMathOperator{\wt}{wt}
\DeclareMathOperator{\coh}{coh}

\renewcommand{\phi}{\varphi}

\DeclareMathOperator{\crk}{crk}

\makeatletter 
\DeclareFontEncoding{LS1}{}{}
\DeclareFontSubstitution{LS1}{stix}{m}{n}
\DeclareMathAlphabet{\mathcaldos}{LS1}{stixscr}{m}{n}
\makeatother

\newcommand{\Module}{\mathbb{M}}

\DeclareMathOperator{\vdim}{vdim}

\DeclareMathOperator{\Perv}{Perv} 
\DeclareMathOperator{\BoMo}{BM} 
\DeclareMathOperator{\Tr}{Tr} 
\DeclareMathOperator{\Tot}{Tot} 
\DeclareMathOperator{\Dol}{Dol} 
\DeclareMathOperator{\Betti}{Betti} 
\DeclareMathOperator{\framed}{fr} 
\DeclareMathOperator{\Higgs}{\mathcal{H}iggs} 
\DeclareMathOperator{\BPS}{\mathtt{BPS}} 
\DeclareMathOperator{\DT}{\mathtt{DT}} 
\DeclareMathOperator{\shBPS}{\mathcal{BPS}} 
\DeclareMathOperator{\JH}{\mathtt{JH}} 
\DeclareMathOperator{\crit}{crit} 
\newcommand{\dd}{\mathbf{d}}

\newcommand{\ff}{\mathbf{f}}
\newcommand{\nn}{\mathbf{n}}
\DeclareMathOperator{\Det}{Det}
\DeclareMathOperator{\poinc}{\mathtt{p}}

\usepackage{mathtools}

\usepackage{hyperref}

\usepackage{algorithmic}
\ifpdf
  \DeclareGraphicsExtensions{.eps,.pdf,.png,.jpg}
\else
  \DeclareGraphicsExtensions{.eps}
\fi

\usepackage{amssymb,stmaryrd}
\DeclareMathOperator{\diag}{diag}

\begin{document}

\newcommand\relatedversion{}
\renewcommand\relatedversion{\thanks{The full version of the paper can be accessed at \protect\url{https://arxiv.org/abs/0000.00000}}} 

\title{BPS cohomology in geometry and representation theory}
    \author{Ben Davison}

\date{}

\maketitle






\begin{abstract} We motivate and survey the theory of BPS invariants of categories and BPS cohomology of stacks, indicating applications in enumerative geometry and representation theory, as well as recent advances.
\end{abstract}

\section{Introduction}
\label{intro_section}
Throughout this paper all schemes and stacks are defined over the complex number field.  Given a variety $X$, its Poincar\'e polynomial $\poinc(X)\coloneqq\sum_{i\in\bZ}(-1)^i\dim(\HO^i(X,\bQ))q^{i/2}$ encodes the ranks of the singular cohomology groups of the underlying complex analytic space of $X$.  This is indeed a polynomial, since $\HO^i(X,\bQ)=0$ for $i\notin [0,\ldots,2\dim(X)]$.

The singular cohomology groups $\HO^i(X,\bQ)$ can be defined by applying the derived global sections functor\footnote{All functors in this paper will be assumed to be derived functors.} $(X\rightarrow \pt)_*$ to the constant sheaf $\bQ_X$ on $X$.  We may also consider more general bounded constructible complexes of sheaves $\cF$ on the space $X$.  In this case, defining $\poinc(X,\cF)\coloneqq\sum_{i\in\bZ}(-1)^i\dim(\HO^i(X,\cF))q^{i/2}$ still yields a (Laurent) polynomial.  An important example of such an $\cF$ is provided by the \textit{dualizing complex} $\bD\bQ_X$.  Since, for a general morphism $f$ of varieties (or stacks) we have a natural isomorphism $f_*\bD\cong \bD f_!$, there is an isomorphism $\HO(X,\bD\bQ_X)\cong \HO_{\mathrm{c}}(X,\bQ)^{\vee}$, where the target is the (graded) dual of the compactly supported cohomology.  So after the substitution $q^{1/2}\mapsto q^{-1/2}$ the polynomial $\poinc(X,\bD\bQ_X)$ is the Poincar\'e polynomial for compactly supported cohomology of $X$.

Now let $\fM$ be an Artin stack.  A motivating example is when $\fM$ is the stack of objects in a category; for concreteness we set $\fM=\fM(\bC[x,y])$ to be the stack of finite-dimensional modules over $\bC[x,y]$.  We define $\fM_n\coloneqq \fM_n(\bC[x,y])$ to be the stack of $n$-dimensional $\bC[x,y]$-modules.  We consider the dualizing complex $\bD\bQ=\bD\bQ_{\fM}$.  Since there is a decomposition $\fM=\coprod_{n\geq 0}\fM_n(\bC[x,y])$, we have $\poinc(\fM,\bD\bQ)=\sum_{n\geq 0}\poinc(\fM_n,\bD\bQ)$.  The right-hand side of this sum does not make sense, since for certain $i$, there are infinitely many $n$ for which the contribution to the $q^i$-coefficient is nonzero.   So consider instead the \textit{partition function}: $\cZ_{\bC[x,y]}(T)=\poinc_{\bN}(\fM,\bD\bQ)\coloneqq \sum_{n\geq 0}\poinc(\fM_n,\bD\bQ)T^n$.  The coefficient of $T^n$ is a formal Laurent power series in $q^{1/2}$.

The partition function $\cZ_{\bC[x,y]}(T)$ can be efficiently expressed in terms of \textit{plethystic exponentials}.  Given a formal power series $f(q,T)=\sum_{n\in \bN,i\in\bZ}f_{i,n}q^{i/2}T^n$ such that $f_{i,n}=0$ for $n=0$, or for $i\ll 0$ depending on $n$, we define $\Exp(f(q,T))\coloneqq\prod_{i,n}(1-q^{i/2}T^n)^{-f_{i,n}}$.  Recursively solving for powers of $T$, we \textit{define} a sequence of formal Laurent power series $\Omega_n(q^{1/2})$, called the \textit{refined BPS invariants} of the category of $\bC[x,y]$-modules, via the equation
\begin{equation}
\label{C2_BPS_def}
\cZ_{\bC[x,y]}(T)=\Exp\left(\sum_{n\geq 1} -\Omega_n(q^{1/2})q^{1/2}(1-q)^{-1}T^n\right).
\end{equation}
Conversely, calculating the formal power series $\Omega_n(q^{1/2})$ for all $n$ determines each of the formal series $\poinc(\fM_n,\bD\bQ)$.  The refined BPS invariants in this example are remarkably simple: we have $\Omega_n(q^{1/2})=-q^{-3/2}$ for all $n$, so that $\cZ_{\bC[x,y]}(T)=\Exp(\sum_{n\geq 1} q^{-1}(1-q)^{-1}T^n)$.  See \cite{MR3032328} for the proof at the level of ``weight'' BPS invariants, and then \cite{MR4661532} for the passage to the version of the result stated here.  In more general problems, the refined BPS invariants are not always so simple to write down, but they are at least always Laurent \textit{polynomials} in $q^{1/2}$.  This is known as \textit{integrality}\footnote{The reason this phenomenon is called integrality will be explained in \S \cref{why_integrality}.  In \cite{MR2851153} this property is also called \textit{admissibility}.} of refined BPS invariants.

In general, for a stack $\fM_{\sA}$ of objects in one of the categories appearing in this paper, the appropriate complex of sheaves $\cF_{\sA}$ to take the derived global sections of depends on what type of category $\sA$ is.  In a very wide class of class of examples beyond $\sA=\bC[x,y]\lmod$, we find that integrality holds for the resulting partition function.  In fact something much stronger holds in our example (and beyond): there is an isomorphism
\begin{align}
\label{C2ic}
\bigoplus_{n\geq 0}\HO(\fM_n(\bC[x,y]),\bD\bQ)\cong \Sym\left(\bigoplus_{n\geq 1} \bQ[2] \otimes\HO(\pt/\bC^*,\bQ)\right),
\end{align}
where $V[j]$ in general denotes the cohomological shift: $V[j]^i=V^{i+j}$.  The left-hand side of \eqref{C2ic} is bigraded: it has a grading by $\bZ_{\coh}\oplus \bN$, where throughout the paper $\bZ_{\coh}$ denotes a copy of $\bZ$ that keeps track of cohomological degrees.  The vector space $V=\bigoplus_{n\geq 1} \bQ[2] \otimes\HO(\pt/\bC^*,\bQ)$ also carries a $\bZ_{\coh}\oplus \bN$-grading: $V_{n}^i=\bQ$ for $(i,n)\in \bZ_{\coh}\oplus \bN$ with $i\in 2\cdot \bZ_{\geq -1}$ and $n\geq 1$, and $V_n^i=0$ otherwise.  The right-hand side of \eqref{C2ic} is  the free commutative $\bZ_{\coh}\oplus \bN$-graded algebra generated by $V$.  It carries a $\bZ_{\coh}\oplus \bN$-grading defined by setting the degree $\lvert a_1\otimes\ldots \otimes a_r\lvert$ of a tensor of homogeneous elements to be $\sum_{s=1}^r\lvert a_s\lvert$.

The relation between \eqref{C2_BPS_def} and \eqref{C2ic} comes from the following equalities: if $V$ is a $\bZ_{\coh}\oplus\bN_{\geq 1}$-graded vector space, with finite-dimensional graded pieces, and such that for each $n\in\bN_{\geq 1}$ we have $V_n^i=0$ for $i\ll0$, then
\begin{equation}
\label{babyExp}
\Exp(\poinc_{\bN}(V))=\poinc_{\bN}(\Sym(V));\quad\quad \Exp(-\poinc_{\bN}(V)q^{1/2}(1-q)^{-1})=\poinc_{\bN}\left(\Sym\left(V\otimes\HO(\pt/\bC^*,\bQ)[-1]\right)\right).
\end{equation}

The \textit{BPS cohomology} in this example is the $\bZ_{\coh}\oplus\bN_{\geq 1}$-graded vector space $\bigoplus_{n\geq 1} \bQ[2]$.  In the first instance, defining the BPS cohomology, calculating it, and proving it fits into an isomorphism as in \eqref{C2ic} solves the integrality problem for the category of $\bC[x,y]$-modules, and gives us a way to express the highly infinite-dimensional vector space $\bigoplus_{n\geq 0}\HO(\fM_n(\bC[x,y]),\bD\bQ)$ in terms of something much more manageable.  As we will see in this paper, the BPS cohomology of various categories also carries interesting algebraic structure, and has deep links with problems in combinatorics, quantum groups, cluster algebras, algebraic geometry and nonabelian Hodge theory.
\section{BPS cohomology of moduli stacks of objects}
\label{SOO}
If $\sC$ is some sufficiently well-behaved dg category, and $\sA^{\sst}$ some subcategory of semistable objects in $\sC$, then we can consider the moduli stack $\fM^{\sst}_{\sA}$ of objects in $\sA^{\sst}$.  As above, we may try to define and construct BPS invariants by looking at the derived global sections of an appropriate sheaf $\cF_{\sA}$ on this stack, taking plethystic exponentials, etc.  Rather than diving directly into the theory for general $\sA^{\sst}$, we start our exposition with a concrete class of examples -- categories of modules over algebras -- before discussing how to generalise to a wider class of categories.  This wider class of categories contains, for example, categories of semistable coherent sheaves on complex projective 3-Calabi--Yau varieties.

Let $A$ be an algebra, and assume that we have fixed a presentation of $A$ of one of the three following forms:
\begin{enumerate}
\item
$A=\bC Q$, the path algebra of a finite quiver over the complex numbers.  This algebra has a $\bC$-basis given by paths in $Q$, including for each vertex $i$ a ``lazy'' path of length zero, beginning and ending at $i$.  Multiplication is given by concatenation of paths, where possible, and zero otherwise.
\item
$A=\bC Q/\langle r_1,\ldots,r_l\rangle$, where each $r_s$ is a linear combination of paths of length at least two in $Q$, such that all the paths in the linear combination $r_s$ have the same starting point and the same endpoint.
\item
$A=\Jac(Q,W)\coloneqq \bC Q/\langle \partial W/\partial a \;\lvert \; a\in Q_1\rangle$, the Jacobi algebra associated to a quiver $Q$ and a potential $W\in \bC Q/[\bC Q,\bC Q]$.  See \cite{ginz} for the definition of the \textit{noncommutative partial derivatives} $\partial W/\partial a$, as well as a wonderful account of the theory of Jacobi algebras and Calabi--Yau algebras in general.
\end{enumerate}
We refer to these presentations as types (1), (2) and (3) from now on.  The cohomology theory we consider, as well as the resulting BPS invariants, will depend on the type of presentation we give $A$.  We denote the vertex set of $Q$ by $Q_0$, and the set of arrows by $Q_1$.  We define $s(a)$ and $t(a)$ to be the source and target of $a\in Q_1$, respectively.  The dimension vector of an $A$-module $\rho$ is the tuple $\dim_{Q_0}(\rho)\coloneqq (\dim(1_i\cdot \rho))_{i\in Q_0}\in \bN^{Q_0}$.  We define the Euler form on $\bN^{Q_0}$
\[
\chi_Q(\dd',\dd'')=\sum_{i\in Q_0}\dd'_i\dd''_i-\sum_{a\in Q_1}\dd'_{s(a)}\dd''_{(t(a)}
\]
and define the symmetrised and antisymmetrised versions
\[
( \dd',\dd'')_Q=\chi_Q(\dd',\dd'')+\chi_Q(\dd'',\dd');\quad\quad \langle \dd',\dd''\rangle_Q=\chi_Q(\dd',\dd'')-\chi_Q(\dd'',\dd').
\]
A quiver $Q'$ is called \textit{symmetric} if $\chi_Q(-,-)$ is symmetric.  In cases (1) and (3) we will initially assume that $Q$ is symmetric.  In case (2), we form a new quiver $Q^+$ by introducing $l$ new arrows $r_s^*$.  The starting point of $r_s^*$ is the same as the endpoint of all of the terms in $r_s$, and the endpoint of $r_s^*$ is the same as the starting point of all of the terms in $r_s$.  We then assume that $Q^+$ is symmetric.

\subsection{Cohomology of stacks of modules}
\label{the_players}

\begin{remark}[Mixed Hodge structures]
\label{MHS_remark}
For each of the three types of algebra above, the cohomology we consider admits a natural lift to the category of \textit{mixed Hodge structures}.  Since we take vanishing cycle cohomology in type (3), we should instead consider \textit{monodromic} mixed Hodge structures and modules as in \cite{MR4132957}.  Where one sees $\bL^{n/2}$ below, which we define to be the $n/2$th tensor power of the pure cohomologically graded mixed Hodge structure $\bL\coloneqq \HO_{\mathrm{c}}(\bA^1,\bQ)$, we should again work with monodromic mixed Hodge structures, since it is in that category that $\bL$ has a tensor square root.  Since, apart from this technicality, so few of the mixed Hodge structures below actually have monodromy\footnote{The essential exceptions are Examples \cref{olqp_example} and \cref{Okkes_example}.}, we will omit the word ``monodromic'' in the hope that doing so reduces confusion.  The paper can be read at the level of cohomologically graded vector spaces, substituting every instance of $\otimes\bL^{n/2}$ for $[-n]$.  Except\footnote{Although, for full disclosure, this failure of purity is one of the most important features of the subject in \S \cref{NAHT_section}!} in \S \cref{NAHT_section} all of the mixed Hodge structures that we encounter will be \textit{pure}, meaning that the refined BPS invariants defined via mixed Hodge structures agree with those defined in a more naive way by simply recording cohomological degree (see \S \cref{BPScotoinv} for the comparison). 
\end{remark}

\begin{enumerate}
\item
For $A=\bC Q$, and $\dd\in\bN^{Q_0}$ a dimension vector, consider the stack $\fM_{\dd}(Q)$ of $\dd$-dimensional $\bC Q$-modules.  It may be written as a global quotient stack $\bA_{\dd}(Q)/\GL_{\dd}$, where
\[
\bA_{\dd}(Q)\coloneqq \prod_{a\in Q_1} \Hom(\bC^{\dd_{s(a)}},\bC^{\dd_{t(a)}});\quad\quad \GL_{\dd}\coloneqq \prod_{i\in Q_0} \GL_{\dd(i)}(\bC)
\]
and $\GL_{\dd}$ acts by simultaneous change of basis.  We define $\cH_{A,\dd}\coloneqq \HO(\fM_{\dd}(Q),\bQ)\otimes\bL^{\chi_Q(\dd,\dd)/2}$.  Since $\fM_{\dd}(Q)$ is a global quotient stack, and $\bA_{\dd}(Q)$ is equivariantly contractible, there is an isomorphism $\cH_{A,\dd}\cong \HO_{\GL_{\dd}}(\pt,\bQ)\otimes\bL^{\chi_Q(\dd,\dd)}$.  The $\bL$-twist is introduced to make the statement of Theorem \cref{int_thm} cleaner.
\item
Fix a presentation $A=\bC Q/\langle r_1,\ldots,r_l\rangle$.  We consider the stack $\fM_{\dd}(A)$ of $\dd$-dimensional $A$-modules.  It is a closed substack of $\fM_{\dd}(Q)$, and is also a global quotient stack.  It is usually singular.  We define $\cH^{\BoMo}_{A,\dd}\coloneqq \HO^{\BoMo}(\fM_{\dd}(A),\bQ)\otimes\bL^{f(\dd)}$, where for an arbitrary space $X$ we define $\HO^{\BoMo}(X,\bQ)=\HO(X,\bD\bQ_X)$, the derived global sections of the Verdier dual of the constant sheaf $\bQ_X$, which if we are pursuing the theory at the level of mixed Hodge structures is upgraded to a mixed Hodge module complex as in \cite{MR1042805}.  Alternatively, we have $\HO^{\BoMo}(X,\bQ)\cong \HO_{\mathrm{c}}(X,\bQ)^{\vee}$, the dual compactly supported cohomology.  We define $f(\dd)=-\chi_Q(\dd,\dd)/2-\sum_{j=1}^l \dd_{s(r_j^*)}\dd_{t(r_j^*)}/2$.  Again, the $\bL$-twists make Theorem \cref{int_thm} cleaner.
\item
Fix a presentation $A=\Jac(Q,W)$.  The potential $W$ defines the function $\Tr(W)$ on $\fM_{\dd}(Q)$, and $\fM_{\dd}(A)=\crit(\Tr(W))$, considered as closed substacks of $\fM_{\dd}(Q)$.  We consider the (perverse exact) vanishing cycle functor $\phi_{\Tr(W)}$, and define $\cH^{\crit}_{A,\dd}\coloneqq \HO(\fM_{\dd}(A),\phi_{\Tr(W)}\bQ_{\fM_{\dd}(Q)})\otimes\bL^{\chi_Q(\dd,\dd)/2}$.
\end{enumerate}
In each type, it is more instructive to consider all of the $\cH_{A,\dd}$ across all dimension vectors simultaneously, as in the following theorem from \cite{MR4132957}.  Given a $\bZ_{\coh}\oplus\bN^{Q_0}$-graded vector space (or mixed Hodge structure) $V=\sum_{\dd\in\bN^{Q_0},n\in \Z}V^n_{\dd}$, we define $\Sym(V)$ to be the free $\bZ_{\coh}\oplus\bN^{Q_0}$-graded \textit{super}commutative algebra generated by $V$: we impose the Koszul sign rule with respect to the cohomological grading, so $v\cdot v'=(-1)^{nn'}v'\cdot v$ if $v,v'$ are of degrees $(n,\dd)$ and $(n',\dd')$ respectively.

\begin{theorem}
\label{int_thm}
Let $A$ be presented as in one of the cases (1), (2) or (3) above, and define $\cH_{A}\coloneqq \bigoplus_{\dd\in\bN^{Q_0}} \cH_{A,\dd}$, one of the three types of $\bZ_{\coh}\oplus \bN^{Q_0}$-graded mixed Hodge structure introduced above.  There is an isomorphism of $\bZ_{\coh}\oplus\bN^{Q_0}$-mixed Hodge structures $\Phi\colon \cH_{A}\cong \Sym\left(\BPS_{A}\otimes \HO(\pt/\bC^*,\bQ)\otimes\bL^{1/2}\right)$ where $\BPS_{A}$ is a $\bZ_{\coh}\oplus\bN^{Q_0}$-graded mixed Hodge structure satisfying the following \textbf{cohomological integrality} property:  For each fixed $\dd\in\bN^{Q_0}$, the $\bZ_{\coh}$-graded mixed Hodge structure $\BPS_{A,\dd}$ is finite-dimensional (after forgetting the $\bZ_{\coh}$-grading).  $\BPS_{A}$ is called the \textbf{BPS cohomology} of the stack of $A$-modules.
\end{theorem}

\subsection{BPS invariants from BPS cohomology}
\label{BPScotoinv}
The definition of plethystic exponentials from \S \cref{intro_section} generalises to multiple variables; we spell out some definitions here, in order to clarify the general relation between BPS cohomology, mixed Hodge structures, BPS invariants, and positivity.  

Fix a rational strongly convex polyhedral cone $\sigma\in \bR^d$, and set $N^+=\sigma\cap \bZ^d$.   Throughout this section, we will set $N^+=\bN^{Q_0}$.  The map $+\colon(N^+)^2\rightarrow N^+$ has finite fibres.  Fix a coefficient ring $R$.  Set $R[T]=R[N^+]$ to be the semigroup algebra of $N^+$ over $R$, which as a $R$-module is freely generated by symbols $T^{\alpha}$ with $\alpha\in N^+$.  Let $I$ be the ideal in $R[T]$ spanned by symbols $T^{\alpha}$ with $\alpha\neq 0$.  Let $R[\![T]\!]$ be the completion of $R[T]$ at $I$, and let $R[\![T]\!]_+\subset R[\![T]\!]$ be the space of formal linear combinations $\sum_{\alpha\in N^+\setminus\{0\}}f_{\alpha}T^{\alpha}$.

We set\footnote{This ring is a somewhat arbitrary choice, chosen to be just flexible to cover everything in this paper.  A more abstract treatment would start by letting $R$ be a general $\lambda$-ring, or the Grothendieck ring of a symmetric tensor category.} $R=\bZ[q_1^{\pm 1/2}\ldots q_{r-1}^{\pm 1/2}](\!(q_r^{1/2})\!)$.  We write elements of $R[\![T]\!]$ as $f=\sum_{\nn\in\bZ^r,\alpha\in N^+}f_{\nn,\alpha}q^{\nn/2} T^{\alpha}$, where we follow the standard notational convention of writing $q^{\nn/2}=\prod_{i=1}^rq_i^{\nn_i/2}$.  Finally, we define 
\[
\Exp\colon R[\![T]\!]_+\rightarrow 1+R[\![T]\!]_+;\quad\quad
\sum_{\nn\in\bZ^r,\alpha\in N^+}f_{\nn,\alpha}q^{\nn/2} T^{\alpha}\mapsto \prod_{\nn\in\bZ^r,\alpha\in N^+}(1-q^{\nn/2}T^{\alpha})^{-f_{\nn,\alpha}}.
\]
This is an isomorphism of groups, where the domain carries the group structure given by addition and the target carries the group structure given by multiplication.

Let $V$ be a $\bZ^{r}\oplus (N^+\setminus\{0\})$-graded vector space, and assume that it also carries a cohomological grading.  Assume that the total dimension of each $\bZ^{r}\oplus N^+$-graded piece is finite-dimensional, and that for fixed $\alpha\in N^+$ we have $V_{(n_1,\ldots,n_r),\alpha}=0$ for $n_r\ll 0$, and for all but finitely many $n_1,\ldots,n_{r-1}$ for fixed $n_r$.  We form the formal Laurent series $\chi_{N^+}(V)=\bigoplus_{(\nn,\alpha)\in\bZ^r\oplus N^+}\sum_{i\in\bZ}(-1)^i\dim(V_{\nn,\alpha}^i)q^{\nn}T^{\alpha}\in R[\![T]\!]_+$.  Then generalising \eqref{babyExp} we have:
\begin{equation}
\label{full_version}
\Exp(\chi_{N^+}(V))=\chi_{N^+}(\Sym(V))
\end{equation}
where $\Sym(V)$ is the $\bZ^r\oplus N^+$-graded free supersymmetric algebra generated by $V$.  We again impose the Koszul sign rule with respect to the cohomological grading.

Given one of the $\bN^{Q_0}$-graded, cohomologically graded vector spaces $\cH_A$ introduced in \S \cref{the_players} we define its \textit{Poincar\'e series}
\[
\poinc_{\bN^{Q_0}}(\cH_A)\coloneqq \sum_{\dd\in\bN^{Q_0};\; i\in \bZ}(-1)^i\dim(\cH_{A,\dd}^i)q^{i/2}T^{\dd}.
\]
Note that the cohomological grading is doing double duty here.  In the general definition of plethystic exponentials, we had a monoid $\bZ^r$ (for now we have $r=1$), which we use to grade vector spaces $V$, and we \textit{also} considered $V$ to have a separate cohomological grading.  We take the alternating sum over cohomological degrees in the definition of $\chi_{N^+}$.  In the definition of Poincar\'e series, however, the $\bZ^r$-grading \textit{is} the cohomological grading.  This yields a cheap \textit{positivity} property: for an arbitrary $\bZ_{\coh}\oplus \bN^{Q_0}$-graded vector space $\cH$ (such that $\cH^i_{\dd}=0$ for $i\ll 0$ depending on $\dd$) we have $\poinc_{\bN^{Q_0}}(\cH)\in \bN(\!(-q^{1/2})\!)[\![T]\!]$.

As in \S \cref{intro_section}, we define \textit{refined BPS invariants} $\Omega^{\coh}_{A,\dd}(q^{1/2})$ by setting
\begin{equation}
\label{coh_BPS_def}
\poinc_{\bN^{Q_0}}(\cH_A)=\Exp\left(\sum_{\dd\in\bN^{Q_0}\setminus\{0\}}-\Omega^{\coh}_{A,\dd}(q^{1/2})q^{1/2}(1-q)^{-1}\right).
\end{equation}
Then \eqref{full_version} along with Theorem \cref{int_thm} implies that $\poinc_{\bN^{Q_0}}(\BPS_A)=\sum_{\dd\in\bN^Q_0/\setminus \{0\}}\Omega^{\coh}_{A,\dd}(q^{1/2})T^{\dd}$.  In particular, all of the refined BPS invariants defined this way are Laurent polynomials in $-q^{1/2}$ with positive coefficients.

As noted in Remark \cref{MHS_remark}, the $\bZ_{\coh}\oplus\bN^{Q_0}$-graded vector space $\cH_A$ can be upgraded to a mixed Hodge structure.  In particular, for each $(i,\dd)$-graded piece $\cH_{A,\dd}^i$ there is an ascending weight filtration $W_{\bullet}\cH_{A,\dd}^i$.  For a mixed Hodge structure $V$ we write $\poinc_{\wt}(V)=\sum_{n\in\bZ}\dim(\Gr_W^n(V))q_1^{n/2}$.  For a $\bZ_{\coh}\oplus\bN^{Q_0}$-graded mixed Hodge structure $V$ we define
\[
\poinc_{\wt,\bN^{Q_0}}(V)\coloneqq\sum_{i\in\bZ;\;\dd\in\bN^{Q_0}}(-1)^i\poinc_{\wt}(V_{\dd}^i)q_2^{i/2}T^{\dd};\quad\quad \chi_{\wt,\bN^{Q_0}}(V)\coloneqq\left(\poinc_{\wt,\bN^{Q_0}}(V)\right)_{\substack{q_1\mapsto q\\q_2\mapsto 1}}.
\]
Then as before, if we define the refined BPS invariants via 
\begin{align}
\poinc_{\wt,\bN^{Q_0}}(\cH_A)&=\Exp\left(\sum_{\dd\in\bN^{Q_0}\setminus\{0\}}-\Omega_{A,\dd}(q_1^{1/2},q_2^{1/2})q^{1/2}(1-q)^{-1}\right)\\
\label{wt_BPS_def}
 \chi_{\wt,\bN^{Q_0}}(\cH_A)&=\Exp\left(\sum_{\dd\in\bN^{Q_0}\setminus\{0\}}-\Omega^{\wt}_{A,\dd}(q^{1/2})q^{1/2}(1-q)^{-1}\right)
\end{align}
where we have set $q=q_1q_2$, we find $\poinc_{\wt,\bN^{Q_0}}(\BPS_A)=\sum_{\dd\in\bN^Q_0\setminus \{0\}}\Omega_{A,\dd}(q_1^{1/2},q_2^{1/2})T^{\dd}$ and $\chi_{\wt,\bN^{Q_0}}(\BPS_A)=\sum_{\dd\in\bN^Q_0 \setminus \{0\}}\Omega^{\wt}_{A,\dd}(q^{1/2})T^{\dd}$.  

Let us assume that $\cH_A$ is \textit{pure}, meaning that $\Gr_W^n(\cH_{A,\dd}^i)=0$ for all $\dd$ if $i\neq n$.  It follows that $\chi_{\wt,\bN^{Q_0}}(\cH_A)=\poinc_{\bN^{Q_0}}(\cH_A)$ and $\chi_{\wt,\bN^{Q_0}}(\BPS_{A})=\poinc_{\bN^{Q_0}}(\BPS_{A})$.  A neat consequence of purity, then, is \textit{positivity} of refined BPS invariants.  Using purity, or an equality between weight and cohomological refinements, to show that polynomials from ``motivic'' DT theory are in fact Poincar\'e polynomials can be used to prove a variety of positivity conjectures, starting with Efimov's proof of the positivity of refined BPS invariants of symmetric quivers \cite{MR2956038}.  We will see a (re)proof of the Kac positivity conjecture via this idea in \S \cref{Kac_section}.  Separate applications of this technique appear in \cite{MR3739230,MR4337972}, to prove positivity conjectures for quantum cluster algebras and quantum theta bases.  

\begin{remark}
It is often easier to calculate the invariants $\Omega^{\wt}_{A,\dd}(q^{1/2})$ than the invariants $\Omega^{\coh}_{A,\dd}(q^{1/2})$.  This is in part because of the motivic nature of weight polynomials (i.e. the fact that for $U\subset X$ an open substack with closed complement $Z$, we have $\chi_{\wt}(\HO^{\BoMo}(U,\bQ))+\chi_{\wt}(\HO^{\BoMo}(Z,\bQ))=\chi_{\wt}(\HO^{\BoMo}(X,\bQ))$.  In many cases (including the example we saw in \S \cref{intro_section}) the way to calculate $\Omega^{\coh}_{A,\dd}(q^{1/2})$ is to first calculate $\Omega^{\wt}_{A,\dd}(q^{1/2})$ and then use purity to show $\Omega^{\coh}_{A,\dd}(q^{1/2})=\Omega^{\wt}_{A,\dd}(q^{1/2})$.  Note that the refined invariants $\Omega^{\wt}_{A,\dd}(q^{1/2})$ agree with the refined invariants considered in \cite{KS1}, while in the absence of purity, $\Omega^{\coh}_{A,\dd}(q^{1/2})$ may not.
\end{remark}
\subsection{BPS sheaves and BPS cohomology}
\label{BPSshco}
Theorem \ref{int_thm} is an existence theorem.  It uniquely determines $\BPS_A$, as a $\bZ_{\coh}\oplus \bN^{Q_0}$-graded mixed Hodge structure, up to isomorphism.  Next we provide a more explicit definition.  Let us start with case (3), so we fix a presentation $A=\Jac(Q,W)$.  Writing $\bC Q=\Jac(Q,0)$, we will simultaneously provide the definition in cases (1) and (3), while the definition for (2) will come in \S \cref{dim_red_sec}.

A key role in the general definition of BPS cohomology is played by \textit{good moduli spaces}.  See \cite{MR3237451} for a general definition of what a good moduli space is, and some of the beautiful ``modern moduli theory'' built around them.  For the global quotient stack $\fM_{\dd}(Q)$, this good moduli space is just the affinization of $\fM_{\dd}(Q)$, the spectrum of the algebra of $\GL_{\dd}$-invariant functions on $\bA_{\dd}(Q)$:
\[
\cM_{\dd}(Q)\coloneqq \Spec(\Gamma(\cO_{\bA_{\dd}(Q)})^{\GL_{\dd}}).
\]
There is a canonical morphism $\JH_{\dd}\colon \fM_{\dd}(Q)\rightarrow \cM_{\dd}(Q)$.  Closed points of $\fM_{\dd}(Q)$ correspond to $\dd$-dimensional $\bC Q$-modules, while closed points of $\cM_{\dd}(Q)$ correspond to semisimple $\dd$-dimensional $\bC Q$-modules.  If $x\in \fM_{\dd}(Q)$ corresponds to the $\dd$-dimensional $\bC Q$-module $\rho$, then $\JH_{\dd}(x)$ corresponds to the semisimplification of $\rho$, that is the direct sum $\bigoplus_{i=1}^l\rho_i/\rho_{i-1}$ where $0=\rho_0\subset \rho_1\subset \ldots\subset \rho_l=\rho$ is a Jordan--H\"older filtration of $\rho$.

The definition of $\BPS_A$ in general rests on the theory of \textit{perverse sheaves}, or their lifts to \textit{mixed Hodge modules} and in particular, generalisations of the Beilinson--Bernstein--Deligne--Gabber decomposition theorem \cite{MR4870047} to good moduli spaces $p\colon \fM\rightarrow \cM$.  The BBDG decomposition theorem states that if $p\colon X\rightarrow Y$ is a projective morphism from a smooth algebraic variety, the direct image $p_*\bQ_X$ can \textit{first} be decomposed as a direct sum of cohomologically shifted perverse sheaves and then \textit{also} that each of these perverse sheaves are semisimple.  Saito's version of the theorem \cite{MR1047415} lifts this statement to mixed Hodge modules.  The following result says that the \textit{first} part\footnote{The second part may actually fail.  That is, the perverse sheaves appearing on the right-hand side of the decomposition in Proposition \cref{vc_decomp_thm} may not be semisimple, and their lifts to mixed Hodge modules may not be pure.} of the decomposition theorem holds for the direct image of the vanishing cycle complex along $\JH_{\dd}$.  See \cite{Kinj_dec} for a generalisation to all good moduli spaces $p\colon\fM\rightarrow \cM$.
\begin{proposition}[\cite{MR4132957}]
\label{vc_decomp_thm}
There is an isomorphism in the derived category of unbounded complexes of perverse sheaves (or their lifts to mixed Hodge modules) on $\cM_{\dd}(Q)$
\[
\JH_{\dd,*}\phi_{\Tr(W)}\bQ\otimes\bL^{\chi_Q(\dd,\dd)/2}\cong \bigoplus_{i\geq 1} \cH^i\left(\JH_{\dd,*}\phi_{\Tr(W)}\bQ\otimes\bL^{\chi_Q(\dd,\dd)/2}\right)[-i].
\]
\label{perv_van_prop}
\end{proposition}
If the isomorphism is being read at the level of constructible complexes, the cohomology functors are the \textit{perverse} cohomology functors of \cite{MR4870047}.  The proposition says that the perverse truncation $\tau_{\leq 0}(\JH_{\dd,*}\phi_{\Tr(W)}\bQ\otimes\bL^{\chi_Q(\dd,\dd)/2})$ vanishes, and $\JH_{\dd,*}\phi_{\Tr(W)}\bQ\otimes\bL^{\chi_Q(\dd,\dd)/2}$ is isomorphic to an (infinite) direct sum of cohomologically shifted mixed Hodge modules (or cohomologically shifted perverse sheaves).
We define the \textit{BPS sheaf}
\begin{equation}
\label{fBPS}
\shBPS_{A,\dd}\coloneqq \cH^0\!\left(\JH_{\dd,*}\phi_{\Tr(W)}\bQ\otimes \bL^{(\chi_Q(\dd,\dd)-1)/2}\right).
\end{equation}
This is, by definition, a mixed Hodge module on $\cM_{\dd}(Q)$, with support contained in $\cM_{\dd}(A)$.  If our algebra $A$ is presented as in case (1) as a path algebra, we simply define $\shBPS_{A,\dd}\coloneqq \cH^0\!\left(\JH_{\dd,*}\bQ\otimes\bL^{(\chi_{Q}(\dd,\dd)-1)/2}\right)$.  

By Proposition \cref{perv_van_prop}, there is an equality $\shBPS_{A,\dd}\otimes\bL^{1/2}=\tau_{\leq 1}(\JH_{\dd,*}\phi_{\Tr(W)}\bQ\otimes\bL^{\chi_{Q}(\dd,\dd)/2})$.  We \textit{define} $\BPS_{A,\dd}\coloneqq\HO(\cM_{\dd}(Q),\shBPS_{A,\dd})$.  We define\footnote{The reason for this notation will become clearer in \S \ref{BPS_Lie_section}.} $\fn^+_{A,\dd}\coloneqq \HO(\cM_{\dd}(A),\tau_{\leq 1}(\JH_{\dd,*}\phi_{\Tr(W)}\bQ\otimes \bL^{\chi_{Q}(\dd,\dd)/2}))=\BPS_{A,\dd}\otimes\bL^{1/2}$.  Implicit in the definition is a presentation of $A$ as a Jacobi/path algebra.  Where we wish to make this presentation explicit, we will write $\BPS_{Q,W},\fn^+_{Q,W}$ etc.  The natural transformation $\tau_{\leq 1}\rightarrow \id$, induces a morphism $\iota \colon \shBPS_{A,\dd}\otimes\bL^{1/2}\rightarrow \JH_{\dd,*}\phi_{\Tr(W)}\bQ\otimes\bL^{\chi_{Q}(\dd,\dd)/2}$.  Applying the derived global sections functor to $\iota$ yields a morphism $\jmath_1\colon \fn^+_{A,\dd}\rightarrow \cH_{A,\dd}$, which is an injection by Proposition \ref{perv_van_prop}.  This defines the BPS cohomology of the category of $A$-modules, along with a canonical embedding (of a cohomological shift of it) as a $\bZ_{\coh}\oplus\bN^{Q_0}$-graded subspace of $\cH_{A}$.
\subsection{Relative integrality}
Having defined the BPS sheaf, we can state a stronger version of the integrality theorem.  First, note that $\cM_A$ is a monoid: points of $\cM_A$ correspond to semisimple $A$ modules, and there is a \textit{finite} \cite[Lemma 2.1]{MR4000572} morphism of schemes $+\colon \cM_A\times \cM_A\rightarrow \cM_A$ taking $(x,y)$ to $x\oplus y$.  This induces a symmetric monoidal structure on $\Perv(\cM_A)$: given a pair of perverse sheaves $\cF,\cG$ the complex $\cF\boxdot\cG\coloneqq +_*(\cF\boxtimes \cG)$ is again perverse, by finiteness of the morphism $+$.  This structure lifts to a symmetric monoidal structure on the category of mixed Hodge modules on $\cM_A$, or complexes of mixed Hodge modules, so we can talk of the free symmetric algebra $\Sym_{\boxdot}(\cF)$ generated by a given object $\cF$ in any one of these categories.  
\begin{theorem}[\cite{MR4132957} Theorem A]
\label{rel_int_thm}
There is an isomorphism in the category of locally bounded below complexes of mixed Hodge modules on $\cM_A$ (where for the domain we take the underlying complex of the algebra object):
\[
\Sym_{\boxdot}\left(\shBPS_{A,\dd}\otimes \HO(\pt/\bC^*,\bQ)\otimes\bL^{1/2}\right)\rightarrow \bigoplus_{\dd\in\bN^{Q_0}}\JH_{\dd,*}\phi_{\Tr(W)}\bQ\otimes\bL^{\chi_Q(\dd,\dd)/2}.
\]
\end{theorem}
Theorem \ref{int_thm} follows from Theorem \ref{rel_int_thm} by passing to derived global sections.
\subsection{Perverse filtration}
\label{perv_filt_section}
By Proposition \ref{vc_decomp_thm} the morphism of complexes $\tau_{\leq n}\left(\JH_{\dd,*}\phi_{\Tr(W)}\bQ\otimes\bL^{\chi_Q(\dd,\dd)/2}\right)\rightarrow \JH_{\dd,*}\phi_{\Tr(W)}\bQ\otimes\bL^{\chi_Q(\dd,\dd)/2}$ has a left inverse for every $n\in\bN$.  Passing to derived global sections, the natural morphism (generalising $\jmath_1$) of cohomologically graded mixed Hodge structures
\[
\fP_n\cH_{A,\dd}\coloneqq \HO\left(\cM_{\dd}(Q),\tau_{\leq n}\JH_{\dd,*}\phi_{\Tr(W)}\bQ\otimes\bL^{\chi_Q(\dd,\dd)/2}\right)\xrightarrow{\jmath_n} \HO\left(\cM_{\dd}(Q),\JH_{\dd,*}\phi_{\Tr(W)}\bQ\otimes\bL^{\chi_Q(\dd,\dd)/2}\right)=\cH_{A,\dd}
\]
is therefore \textit{injective}.  We call $\fP_n\cH_{A,\dd}$ the $n$th piece of the perverse filtration of $\cH_{A,\dd}$.  Note that $\fn^+_A\coloneqq\fP_1\cH_A$.

\subsection{BPS cohomology of stacks of objects in 3-Calabi--Yau categories}
\label{general_BPS}
Let $\sC$ be a 3-Calabi--Yau (3CY) dg category.  We refer to \cite{MR3911626} for exact definitions, merely recalling that a consequence of the definition is that for objects $x,y$ of $\sC$ satisfying a local compactness condition\footnote{Local compactness of $z$ is the condition that $\Hom_{\sC}(z',z)$ is finite-dimensional for all $z'$.  In the category of coherent sheaves on a noncompact 2CY surface, for example, local compactness of $\cG$ is ensured by the condition that the support of $\cG$ be projective.} there is a bifunctorial pairing $\Ext^i(x,y)\cong \Ext^{3-i}(y,x)^{\vee}$, in analogy with Serre duality for coherent sheaves on a 3-Calabi--Yau variety.  We pick an Abelian heart $\sA\subset \sC$, and we let $\sA^{\sst}\subset \sA$ be a full Abelian subcategory of semistable objects, such that $\fM_{\sA}^{\sst}\subset \fM_{\sC}$, the open substack of objects contained in $\sA^{\sst}$, has a good moduli space $p\colon \fM_{\sA}^{\sst}\rightarrow \cM_{\sA}^{\sst}$ in the sense of \cite{MR3237451}, and the Euler form $\chi_{\sC}(x,y)=\sum_{i\in \bN}(-1)^{i}\dim(\Ext^i(x,y))$ is symmetric on objects $x,y$ of $\sA^{\sst}$.  Rather than expanding upon these definitions, we mention our motivating examples:
\begin{enumerate}
\item
$\sA=\Jac(Q,W)\lmod$ for $Q$ a quiver.  If $Q$ is symmetric, we may pick $\sA^{\sst}=\sA$.  Otherwise we pick a generic stability condition and let $\sA^{\sst}$ be the category of semistable modules of fixed slope as defined in \S \cref{CWC}.
\item
$\sA=\Coh(X)$ for $X$ a quasiprojective 3CY variety.  We pick a generic polarization $H$ on a compactification of $X$, and fix $\sA^{\sst}$ to be the category of $H$-semistable sheaves on $\overline{X}$ with support in $X$, with fixed normalised Hilbert polynomial with respect to $H$.
\item
$\sA^{\sst}=\sA=\bC[\pi_1(M)]\lmod$ for $M$ an orientable closed 3-manifold without boundary.
\end{enumerate}
In each of the above examples, the stack of objects is locally modelled as the critical locus of a function on a smooth stack \cite{MR3352237}.  Our discussion of $\fM(\Jac(Q,W))$ thus provides a \textit{local model} for the global geometry of the stack $\fM_{\sA}^{\sst}$.  In particular, on an atlas of charts $Z_i$ for the stack $\fM_{\sA}^{\sst}$ we can define the vanishing cycle sheaf $\phi_{f_i}\bQ_{U_i}\otimes\bL^{-\dim(U_i)/2}$, where $Z_i=\crit(f_i)$, and $f_i$ is a function on a smooth stack $U_i$.  In all of the above examples, it is possible \cite{MR4224043,NS23} to provide the stack of objects $\fM_{\sA}^{\sst}$ with an extra structure called \textit{orientation data}, which provides the required data to glue these vanishing cycle sheaves across all of the critical charts, defining the mixed Hodge module $\phi_{\DT}$ on $\fM_{\sA}^{\sst}$ -- see \cite{MR3352237} again.

In case (1), the morphism $\JH$ \textit{is} the good moduli space morphism $p$ (at least if $Q$ is symmetric), and $\phi_{\DT}=\phi_{\Tr(W)}\bQ\otimes\bL^{-\dim(\fM(Q))/2}$, while in case (2) $\JH$ models $p$ analytically locally \cite{MR3811778}.  Note that perverse truncation functors commute with taking open restrictions.  It is thus natural to define the BPS sheaf and BPS cohomology, generalising \eqref{fBPS}, as in \cite{MR4564413}:
\[
\shBPS^{\sst}_{\sA}\coloneqq \cH^0(p_*\phi_{\DT}\otimes\bL^{-1/2});\quad\quad \BPS^{\sst}_{\sA}\coloneqq \HO(\cM_{\sA}^{\sst},\shBPS^{\sst}_{\sA}).
\]
Of course, these definitions are only justified if we can prove that $\BPS^{\sst}_{\sA}$ and $\shBPS^{\sst}_{\sA}$ satisfy the analogue of Theorems \ref{int_thm} and \ref{rel_int_thm}, i.e. the cohomological integrality theorem.  Under certain technical conditions on the orientation data and stack $\fM_{\sA}$, this version of the integrality theorem is proved in \cite{BDINKP}, i.e. it is shown that there is an isomorphism of mixed Hodge module complexes
\begin{equation}
\label{genint}
\Sym_{\boxdot}\left(\shBPS^{\sst}_{\sA}\otimes\HO(\pt/\bC^*,\bQ)\otimes\bL^{1/2}\right)\xrightarrow{\cong} p_*\phi_{\DT}
\end{equation}
where the morphism uses \textit{cohomological Hall induction}, discussed in \S \ref{CoHI}.  See \cite{BDINKP} for full details.
\subsection{Dimensional reduction}
\label{dim_red_sec}
Given $A$ with presentation of type (2), define $B=\Jac(Q^+,W)$, where $W=\sum_{j=1}^l r_jr_j^*$.  Let $\pi\colon \fM_{\dd}(Q^+)\rightarrow \fM_{\dd}(Q)$ be the forgetful morphism taking $\bC Q^+$-modules to their underlying $\bC Q$-modules.  By the \textit{dimensional reduction isomorphism} \cite[Appendix A]{MR3667216} there is a natural isomorphism $\pi_*\phi_{\Tr(W)}\bQ\otimes\bL^{\chi_{Q^+}(\dd,\dd)/2}\cong(\bD\bQ_{\fM_{\dd}(A)})\otimes\bL^{f(\dd)}$.  Passing to derived global sections, we find that there is a natural isomorphism $\Psi\colon \cH^{\BoMo}_A\cong \cH^{\crit}_B$.  We define $\BPS_A\coloneqq \BPS_B\coloneqq \BPS_{Q^+,W}$, and we define the embedding $\BPS_A\subset \cH_A$ by composing $\BPS_B\hookrightarrow \cH^{\crit}_B\xrightarrow{\Psi^{-1}}\cH^{\BoMo}_A$.  By Theorem \ref{int_thm} applied to $B$, there is an isomorphism
\[
\cH^{\crit}_{B}\cong \Sym\left(\BPS_{B}\otimes \HO(\pt/\bC^*,\bQ)\otimes\bL^{1/2}\right),
\]
so that via $\Psi^{-1}$ we deduce Theorem \ref{int_thm} also for $\cH^{\BoMo}_A$.

Let $S$ be a smooth quasiprojective surface, let $X=\Tot_S(\omega_S)$ be the total space of the canonical bundle on $S$, and denote by $\varpi\colon X\rightarrow S$ the projection.  Then $X$ is a noncompact 3CY variety.  In this context, a more sophisticated form of the dimensional reduction isomorphism due to Kinjo \cite{MR4379049} applies: if $q\colon \fM_{\Coh(X)}\rightarrow \fM_{\Coh(S)}$ is the morphism of moduli stacks of objects induced by $\varpi_*$, then $q_*\phi_{\DT}\cong (\bD\bQ_{\fM_{\Coh(S)}})\otimes \bL^{\vdim/2}$, where $\vdim$ at a point representing a coherent sheaf $\cF$ is given by $-\chi(\cF,\cF)$.  The sheaf $\phi_{\DT}$ is constructed as in \S \ref{general_BPS}.  If $S$ is moreover a 2CY surface\footnote{We assume in this case that the hyperplane class $H$ is pulled back from $\overline{S}$.}, then $q^{-1}(\fM^{\sst}_{\Coh(S)})=\fM^{\sst}_{\Coh(X)}$, so that after passing to derived global sections we obtain an isomorphism $\HO(\fM^{\sst}_{\Coh(X)},\phi_{\DT})\cong \HO^{\BoMo}(\fM^{\sst}_{\Coh(S)},\bQ)\otimes \bL^{\vdim/2}\eqcolon \cH^{\sst}_{\Coh(S)}$.  As above, we use this isomorphism to define the BPS cohomology of the stack of semistable coherent sheaves on $S$ by transporting the BPS cohomology of the stack of semistable coherent sheaves on $X$, and deduce cohomological integrality for $\cH^{\sst}_{\Coh(S)}$.

\section{Calculating refined BPS invariants}
\subsection{Examples from flopping curves}
It is sometimes possible to compute BPS sheaves by hand, using the following theorem of Meinhardt and Reineke\footnote{In fact they proved the theorem at the level of classes in the Grothendieck ring of mixed Hodge modlues on $\cM_{\dd}(Q)$.  The version stated here is then a straightforward consequence of \textit{purity}, see \cite{MR4132957}.}:
\begin{theorem}[\cite{MR4000572}]
\label{MR_thm1}
Let $Q$ be a symmetric quiver, and let $\dd\in\bN^{Q_0}$.  There is an isomorphism
\[
\shBPS_{Q,\dd}\coloneqq \begin{cases} \IC_{\cM_{\dd}(Q)}&\textrm{if there exists a simple }\dd\textrm{-dimensional }\bC Q\textrm{-module}\\
0&\textrm{otherwise.}
\end{cases}
\]
\end{theorem}
Here, $\IC_{\cM_{\dd}(Q)}$ is the intersection complex: the unique simple mixed Hodge module (or perverse sheaf) extending the constant mixed Hodge module $\bQ_{\cM^{\mathrm{simp}}_{\dd}(Q)}\otimes \bL^{\chi_{Q}(\dd,\dd)/2}$ on the (smooth) locus of simple $\dd$-dimensional $\bC Q$-modules.  Via the natural isomorphism\footnote{See \cite[Section 4.3]{MR4000572} for this isomorphism, and \cite[Appendix A]{Kinj_dec} for the generalisation to \textit{all} good moduli space morphisms.} $\JH_{\dd,*}\phi_{\Tr(W)}\cong \phi_{\Tr(W)}\JH_{\dd,*}$ we deduce that $\shBPS_{Q,W,\dd}$ is isomorphic to $\phi_{\Tr(W)}\IC_{\cM_{\dd}(Q)}$ if there exists a simple $\dd$-dimensional $\bC Q$-module, and is zero otherwise.  
\begin{example}
\label{qdlog_example}
Let $Q^{(0)}$ be the quiver consisting of one vertex and no arrows.  Then $\bC Q^{(0)}=\bC$, and a $\bC Q^{(0)}$-module is just a $\bC$-vector space.  There is a simple $d$-dimensional $\bC Q^{(0)}$-module if and only if $d=1$.  Then $\cM_d(Q^{(0)})=\pt$ and $\IC_{\cM_{\dd}(Q^{(0)})}=\bQ_{\pt}$, so 
\[
\BPS_{A,d}=\begin{cases} \bQ & \textrm{if }d=1\\
0& \textrm{otherwise.}
\end{cases}
\]
We deduce that there is an isomorphism $\cH_Q\cong \Sym\left(\HO(\pt/\bC^*,\bQ)\otimes\bL^{1/2}\right)$.  Passing to underlying $\bZ_{\coh}\oplus \bN$-graded vector spaces, the right-hand side is the underlying $\bZ_{\coh}\oplus\bN$-graded vector space of a free supercommutative algebra generated by a vector space of dimension $1$ in degrees $(1,2n+1)$ for $n\in\bN$, and dimension $0$ in all other degrees.  This is a free \textit{exterior} algebra.  Passing to generating functions, we find the \textit{quantum dilogarithm}
\[
\poinc_{\bN}(\cH_{Q^{(0)}})=\sum_{n\geq 0} (-q^{1/2})^{n^2}\prod_{j=1}^n(1-q^{2j})^{-1}T^n=\prod_{n\in \bN}(1-Tq^{n+1/2})\eqqcolon\bE(T).
\]
More generally we define $\bE(f(q^{1/2},T))\coloneqq \Exp(f(q^{1/2},T)(-q^{1/2}(1-q)^{-1}))$.
\end{example}

\begin{example}
\label{boson_example}
Let $Q^{(1)}$ be the quiver consisting of one vertex and one loop $x$.  Then $\bC Q^{(1)}=\bC[x]$.  The only simple $\bC Q^{(1)}$-modules are one-dimensional.  Moreover $\cM_1(Q^{(1)})\cong\bA^1$ and $\IC_{\cM_1(Q^{(1)})}=\bQ_{\bA^1}\otimes\bL^{-1/2}$.  We deduce that there is an isomorphism $\cH_{Q^{(1)}}\cong \Sym(\HO(\pt/\bC^*,\bQ))$.  This is a free symmetric algebra (it is concentrated in even cohomological degrees) with generating function $\bE(-q^{-1/2}T)$.
\end{example}

\begin{example}
\label{olqp_example}
Let $Q^{(1)}$ be the quiver from the previous example.  Let $W=x^{d+1}$.  Then $\Tr(W)$ induces the function $y^{d+1}$ on the space $\cM_1(Q^{(1)})\cong \bA^1=\Spec(\bC[y])$.  The vanishing cycle cohomology carries a nontrivial monodromy action -- we consider this problem just at the level of vector spaces.  There is an isomorphism $\phi_{y^{d+1}}\bQ_{\bA^1}[1]\cong\bQ_0^{\oplus d}$, and Theorem \ref{MR_thm1} implies that there are isomorphisms $\cH_{Q^{(1)},x^{d+1}}\cong \Sym\left((\HO(\pt/\bC^*,\bQ)[-1])^{\oplus d}\right)\cong \cH_{Q^{(0)}}^{\otimes d}$.
\end{example}

\begin{example}
\label{Okkes_example}
This example is \textit{much} more intricate, and is worked out by van Garderen in \cite{vGard}.  Let $Q^{(2)}$ be the quiver with one vertex and two loops, labelled $x,y$.  Pick $a\in\bN_{\geq 2}$ and $b\in \bN_{\geq 1}\cup\{\infty\}$.  Let $W=x^2y-f_{a,b}(y)$, where 
\[
f_{a,b}(y)=\begin{cases} y^{2a}&b=\infty\\
y^{2a}+y^{2b+1}&b\neq \infty.\end{cases}
\]
Set $A=\Jac(Q^{(2)},W)$.  Again, there is nontrivial monodromy when we calculate the vanishing cycles functor.  It is proved in \cite{vGard} that we have
\[
\BPS_{A,d}=\begin{cases} \bQ^{\oplus \min\{2a+1,2b+2\}} & d=1\\
\bQ^{\oplus a-1}& d=2\\
0& d=3\end{cases}
\]
and the BPS cohomology is pure.  In particular, the (refined) BPS invariants for $A$ are $(\min\{2a+1,2b+2\},a-1,0,0,\ldots)$.
\end{example}
Despite their purely algebraic appearance, all four of the above examples come from 3-dimensional geometry, more precisely the contribution of a rational curve $C$ in different Calabi--Yau threefolds to the BPS invariants $\Omega^{\wt}_{n[C]}(q^{1/2})$ of those threefolds: see \cite{MR2420017,MR3414491,MR4819391,MR3504176} and references therein for some of this very rich story.  One consequence of Example \ref{Okkes_example} is that, even quite deep into this story, the refined BPS invariants are constant, so it is natural to make the following:
\begin{conjecture}
Let $A$ be the contraction algebra of a rational flopping curve $C$: a certain type of Jacobi algebra defined in \cite{MR3504176}.  Then for all $d\in\bN$, the refined BPS invariants $\Omega_{A,d}^{\coh}(q^{1/2})=\Omega_{A,d}^{\wt}(q^{1/2})$ are constant, and are equal to the $d$th Gopakumar--Vafa invariant of $C$.
\end{conjecture}

\subsection{Kac polynomials as refined BPS invariants}
\label{Kac_section}
In \cite{MR718127} Victor Kac proved the following beautiful result.  See \cite{MR3966814} for a more recent survey of Kac polynomials and their generalisations.
\begin{theorem}
Let $Q$ be an arbitrary quiver (finite, but not necessarily symmetric) and pick $\dd\in\bN^{Q_0}$.  Then there exists a polynomial $\kac_{Q,\dd}(t)\in\bZ[t]$ such that if $\bF$ is a finite field with $q$ elements, the number of isomorphism classes of absolutely indecomposable\footnote{A $\bF Q$-module $\rho$ is called absolutely indecomposable if $\rho\otimes_{\bF}\overline{\bF}$ is indecomposable as a $\overline{\bF} Q$-module.} $\dd$-dimensional $\bF Q$-modules is equal to $\kac_{Q,\dd}(q)$.
\end{theorem}
Given a finite quiver $Q$, we form the \textit{doubled} quiver $\overline{Q}$ by setting $\overline{Q}_0=Q_0$ and $\overline{Q}_1=Q_1\coprod Q_1^{\opp}$, where $Q_1^{\opp}$ contains an arrow $a^*$ for every arrow $a\in Q_1$, and we give $a^*$ the opposite orientation to $a$.  We form the \textit{tripled} quiver $\tilde{Q}$ from $\overline{Q}$ by adding a loop $\omega_i$ for every $i\in Q_0$, with $s(\omega_i)=t(\omega_i)=i$.  We consider the \textit{preprojective algebra} $\Pi_Q\coloneqq \bC \overline{Q}/\langle \sum_{a\in Q_1} [a,a^*]\rangle$.  The \textit{canonical cubic potential} for $\tilde{Q}$ is given by $\tilde{W}=(\sum_{a\in Q_1}[a,a^*])(\sum_{i\in Q_0}\omega_i)$.  If we set $A=\Jac(\tilde{Q},\tilde{W})$, there is an isomorphism $A\cong \Pi_Q[\omega]$ sending $\sum_{i\in Q_0}\omega_i$ to $\omega$.

The following is a combination of \cite{Moz11}, which calculates the motivic BPS invariants and hence $\Omega^{\wt}_{A,\dd}(q^{1/2})$, and \cite{MR4661532}, in which purity of the cohomological BPS invariants is proved.
\begin{theorem}
\label{Kac_to_BPS}
For an arbitrary quiver $Q$ and dimension vector $\dd\in\bN^{Q_0}$ we have equalities
\[
\kac_{Q,\dd}(q^{-1})=-q^{1/2} \Omega^{\wt}_{A,\dd}(q^{1/2})=-q^{1/2}\Omega^{\coh}_{A,\dd}(q^{1/2}).
\]
\end{theorem}
Since $\kac_{Q,\dd}(q^{-1})$ is a polynomial in $q^{-1}$, it follows that $\fn^+_{A}$ is supported entirely in even negative cohomological degrees.  Since $-q^{1/2}\Omega^{\coh}_{A,\dd}(q^{1/2})$ has only positive coefficients, it follows that the coefficients of $\kac_{Q,\dd}(q^{-1})$ are positive.    This positivity was originally conjectured by Kac, and then originally proved by Hausel, Letellier and Rodriguez-Villegas in \cite{MR3034296}.

\subsection{Cohomological wall-crossing and quantum dilogarithm identities}
\label{CWC}
We discuss algebras of type (3), and relax the condition that $Q$ is symmetric.  Setting $W=0$, this includes as a special case those algebras given a presentation of type (1).

We introduce a \textit{stability condition} $\zeta\in\bQ^{Q_0}$.  We define the \textit{slope} of a nonzero $A$-module $\rho$ to be
\[
\mu^{\zeta}(\rho)\coloneqq \frac{\zeta\cdot \dim_{Q_0}(\rho)}{\dim(\rho)}.
\]
A nonzero $A$-module $\rho$ is called $\zeta$-\textit{semistable} if for all $0\neq \rho'\subset \rho$ we have $\mu^{\zeta}(\rho')\leq \mu^{\zeta}(\rho)$.  We denote by $\fM_{\dd}^{\zeta\sst}(A)$ the stack of $\dd$-dimensional $\zeta$-semistable $A$-modules.  It is equal to the substack $\crit(\Tr(W))$, where $\Tr(W)$ is considered as a function on the open substack $\fM_{\dd}^{\zeta\sst}(Q)\subset \fM_{\dd}(Q)$ of $\zeta$-semistable $\bC Q$-modules.  We consider $\varphi=\phi_{\Tr(W)}\bQ_{\fM_{\dd}^{\zeta\sst}(Q)}\otimes\bL^{\chi_Q(\dd,\dd)/2}$, the mixed Hodge module of vanishing cycles, as a mixed Hodge module (or just perverse sheaf) on $\fM_{\dd}^{\zeta\sst}(A)$.  We set $\cH_{A,\dd}^{\crit,\zeta}\coloneqq \HO(\fM_{\dd}^{\zeta\sst}(A),\varphi)$.  For fixed $\theta\in\bQ$, set $\Lambda_{\theta}^{\zeta}\coloneqq \{\dd\in\bN^{Q_0}\setminus \{0\}\;\lvert\; \mu^{\zeta}(\dd)=\theta\}\cup\{0\}$.  We define $\cH^{\crit,\zeta}_{A,\theta}=\bigoplus_{\dd\in\Lambda^{\zeta}_{\theta}}\cH_{A,\dd}^{\crit,\zeta}$.  On the category of $\bN^{Q_0}$-graded mixed Hodge structures, we consider the following twisted monoidal structure as in \cite[Section 3.2]{MR4132957}
\[
\cF\otimes^{\twisted}\cG\coloneqq \bigoplus_{\dd\in\bN^{Q_0}}\left(\bigoplus_{\dd'+\dd''=\dd}\bL^{\langle\dd'',\dd'\rangle/2}\otimes \cF_{\dd'}\otimes \cG_{\dd''}\right).
\]
\begin{theorem}[\cite{MR4132957} Theorem B]
\label{wall_crossing_theorem}
There is an isomorphism of $\bZ_{\coh}\oplus\bN^{Q_0}$-graded mixed Hodge structures
\begin{equation}
\label{HNfun}
\cH^{\crit}_A\cong \bigotimes_{\infty \xrightarrow{\theta} -\infty}\cH^{\crit,\zeta}_{A,\theta}.
\end{equation}
It follows that if $\zeta'$ is another stability condition, we have an isomorphism $\bigotimes_{\infty \xrightarrow{\theta} -\infty}\cH^{\crit,\zeta}_{A,\theta}\cong \bigotimes_{\infty \xrightarrow{\theta} -\infty}\cH^{\crit,\zeta'}_{A,\theta}$.
\end{theorem}
Due to this independence as we cross walls in the space of stability conditions we call this the \textit{cohomological wall-crossing theorem}.
\begin{example}
\label{A2example}
Let $Q$ be the $A_2$ quiver: it has two vertices $0$ and $1$, along with a single arrow $a$ with $s(a)=0$ and $t(a)=1$.  Set $A=\bC Q$.  Pick $\zeta^+\in\bQ^{Q_0}$ with $\zeta^+_0=0$ and $\zeta^+_1=1$.  The only $\zeta^+$-semistable modules are either entirely supported at $0$ or at $1$.  As in Example \ref{qdlog_example}, we find $\cH^{\zeta^+}_{A,\theta}\cong  \Sym\left( V_{\theta}\otimes\HO(\pt/\bC^*,\bQ)\otimes \bL^{1/2}\right)$ with $V_{\theta}=\bQ_{(1,0)}$ a copy of $\bQ$ in $\bN^{Q_0}$-degree $(1,0)$ if $\theta=0$, a copy $\bQ_{(0,1)}$ of $\bQ$ in $\bN^{Q_0}$-degree $(0,1)$ if $\theta=1$, and $V_{\theta}=0$ otherwise.

Now consider $\zeta^-=-\zeta^+$.  There are strictly more $\zeta^-$-semistable modules.  Direct sums of $\bC Q$, considered as a $\bC Q$-module, are $\zeta^-$-stable modules of slope $-1/2$.  We find $\cH^{\zeta^-}_{A,1/2}\cong  \Sym\left( \bQ_{(1,1)}\otimes\HO(\pt/\bC^*,\bQ)\otimes \bL^{1/2}\right)$.  The wall-crossing isomorphism from Theorem \ref{wall_crossing_theorem} reads
\begin{align*}
&\Sym\left( \bQ_{(0,1)}\otimes^{\twisted}\HO(\pt/\bC^*,\bQ)\otimes \bL^{1/2}\right)\otimes^{\twisted} \Sym\left( \bQ_{(1,0)}\otimes\HO(\pt/\bC^*,\bQ)\otimes \bL^{1/2}\right)\cong\\ &\Sym\left( \bQ_{(1,0)}\otimes\HO(\pt/\bC^*,\bQ)\otimes \bL^{1/2}\right)\otimes^{\twisted} \Sym\left( \bQ_{(1,1)}\otimes\HO(\pt/\bC^*,\bQ)\otimes \bL^{1/2}\right)\otimes^{\twisted}\Sym\left( \bQ_{(0,1)}\otimes\HO(\pt/\bC^*,\bQ)\otimes \bL^{1/2}\right).
\end{align*}
Passing to generating series, this becomes the \textit{pentagon identity} $\bE(T^{(0,1)})\bE(T^{(1,0)})=\bE(T^{(1,0)})\bE(T^{(1,1)})\bE(T^{(0,1)})$, relating quantum dilogarithms (as introduced in Example \ref{qdlog_example}).  Here, the formal variables $T^{\dd}$ satisfy the relation $T^{\dd'}\cdot T^{\dd''}=(-q^{1/2})^{\langle \dd'',\dd'\rangle_Q}T^{\dd'+\dd''}$, and thus only $q$-commute, due to the non-symmetry of $Q$.
\end{example}

\begin{example}
\label{theta_example}
Pick numbers $l_0,l_1\in\{0,1\}$ and $k\in\bZ_{\geq 1}$. Let $Q$ be the quiver with two vertices $0,1$, with $l_i$ loops at vertex $i$, and $k$ arrows from $0$ to $1$.  Taking characteristic functions of the right-hand side of \eqref{HNfun} for the stability condition $\zeta^+$ from Example \ref{A2example} yields $\bE((-q^{1/2})^{-l_1}T^{(0,1)})\bE((-q^{1/2})^{-l_0}T^{(1,0)})$, the product of two generalised quantum dilogarithms.  In contrast with Example \ref{A2example}, for general choices of $l_0,l_1,k$ the corresponding factorisation for $\zeta^-$ has infinitely many nontrivial terms $\bE(f_{\theta})$ for $\theta\in\bQ$ and $f_{\theta}=\sum_{d_1/d_2=\theta}f_{d_1,d_2}(q^{1/2})T^{(d_1,d_2)}$.  One can show, by identifying coefficients of $f_{\theta}$ with Poincar\'e polynomials of refined BPS invariants, that each $f_{d_1,d_2}(q^{1/2})$ is in $\bN[-q^{1/2}]$, and using recursive scattering diagram techniques this is enough to prove the positivity of structure constants for rings of quantum theta functions: see \cite{MR4337972}.  
\end{example}
A curious generalisation of Example \ref{theta_example} also appears to hold, and using the same scattering diagram techniques this should entail a \textit{Lefschetz type} property for quantum theta functions: see \cite{MR4337972}.  It remains a conjecture since the statement has so far resisted effective translation into the language of BPS cohomology.  For $n\in\bZ_{\geq 1}$ define the (signed) \textit{quantum numbers} $[\pm n]_q=(-1)^{n+1}(q^{n/2}-q^{-n/2})/(q^{1/2}-q^{-1/2})$.
\begin{conjecture}[\cite{MR4337972} Conjecture 1.4]
Pick $n_1,n_2,k\in\bZ_{\geq 1}$.  Write $\bE([\pm n_1]_qT^{(0,1)})\bE([\pm n_2]_qT^{(1,0)})=\prod_{\infty\xrightarrow{\theta}-\infty}\bE(g_{\theta})$ with $g_{\theta}=\sum_{d_1/d_2=\theta}g_{d_1,d_2}(q^{1/2})T^{(d_1,d_2)}$ where $T^{\dd}T^{\dd'}=(-q^{1/2})^{-\langle \dd',\dd''\rangle_Q}T^{\dd+\dd'}$ with $Q$ as in Example \ref{theta_example} (set $l_0=l_1=0$).  Then each $g_{d_1,d_2}(q^{1/2})$ is a sum of elements of $\{[\pm n]_q\;\lvert\;n\geq 1\}$.
\end{conjecture}

\subsection{Why ``integrality''?}
\label{why_integrality}
Let $Q$ be a quiver, which we assume to be symmetric.  Pick a dimension vector $\ff\in\bN^{Q_0}\setminus\{0\}$.  We form the \textit{framed quiver} $Q_{\ff}$ as follows: we set $(Q_{\ff})_0=Q_0\cup\{\infty\}$, and set $(Q_{\ff})_1=Q_1\cup\{r_{i,n}\;\lvert\; i\in Q_0, 1\leq n\leq \ff_i\}$.  We set $s(r_{i,n})=\infty$ and $t(r_{i,n})=i$.  Given a potential $W$ for $Q$, i.e. a linear combination of cyclic paths in $Q$, the same linear combination of cyclic paths defines a potential for $Q_{\ff}$, which we continue to denote $W$.  We pick the stability condition $\zeta$ with $\zeta_i=0$ for all $i\in Q_0$ and $\zeta_{\infty}=1$.  Given a dimension vector $\dd\in\bN^{Q_0}$ and $n\in\bN$ we define $(\dd,n)$, extending $\dd$ by setting $(\dd,n)_{\infty}=n$.  A $(\dd,1)$-dimensional $\bC Q_{\ff}$-module $\rho$ is $\zeta$-semistable if and only if it is generated, as a $\bC Q_{\ff}$-module, by the one-dimensional vector space $1_{\infty}\cdot \rho$.  There is an isomorphism of stacks 
\[
\bA_{(\dd,1)}^{\zeta\sst}(Q_{\ff})/\GL_{(\dd,1)}\cong (\bA_{(\dd,1)}^{\zeta\sst}(Q_{\ff})/\GL_{\dd})/\bC^*
\]
where $\GL_{\dd}$ acts by change of basis on each of the vector spaces $\bC^{\dd_i}$ for $i\in Q_0$, and $\bC^*$ is embedded as the diagonal in $\GL_{(\dd,1)}$ and thus acts trivially.  The quotient $(\bA_{(\dd,1)}^{\zeta\sst}(Q_{\ff})/\GL_{\dd})$ is in fact a \textit{smooth scheme}, which we denote $\cM^{\ff\framed}_{\dd}(Q)$.  Taking vanishing cycle cohomology this gives an isomorphism
\[
\cH^{\zeta}_{Q_{\ff},W,(\dd,1)}\cong  \Module^{\zeta}_{Q,W}(\ff,\dd)\otimes\bL^{1/2}\otimes \HO(\pt/\bC^*,\bQ);\quad\quad \Module^{\zeta}_{Q,W}(\ff,\dd)\coloneqq\HO(\cM^{\ff\framed}_{\dd}(Q),\phi_{\Tr(W)}\bQ\otimes \bL^{-d/2})
\]
where $d=\dim(\cM^{\ff\framed}_{\dd}(Q))$.  It follows that $\chi_{\wt}(\cH^{\zeta}_{Q_{\ff},W,(\dd,1)})=-\chi_{\wt}(\Module^{\zeta}_{Q,W}(\ff,\dd))q^{1/2}(1-q)^{-1}$.

Now consider instead the stability condition $\zeta'=-\zeta$.  A $\bC Q_{\ff}$-module is $\zeta'$-semistable if and only if it is either entirely supported at $Q$, or entirely supported at $\infty$.  We consider the wall-crossing isomorphism from Theorem \ref{wall_crossing_theorem}, and restrict it to those summands corresponding to dimension vectors $(\dd,1)\in\bN^{(Q_{\ff})_0}$:
\[
\cH_{Q,W}\otimes^{\twisted} \bigoplus_{\dd\in\bN^{Q_0}}\cH^{\zeta}_{Q_{\ff},W,(\dd,1)}\cong \cH_{Q_{\ff},W,(0\ldots,0,1)}\otimes^{\twisted}\cH_{Q,W}
\]
where we have identified $\cH_{Q,W,\dd}$ with $\cH_{Q_{\ff},W,(\dd,0)}$.  Note that $\fM_{(0,\ldots,0,1)}(Q_{\ff})\cong\pt/\bC^*$.  Passing to generating functions and dividing by $-q^{1/2}(1-q)^{-1}$ we calculate
\[
\chi_{\wt,\bN^{Q_0}}(\overline{\Module}^{\zeta}_{Q,W}(\ff,-))= T_{\infty}^{-1}(\chi_{\wt,\bN^{Q_0}}(\cH_{Q,W})^{-1}\cdot T_{\infty}\cdot \chi_{\wt,\bN^{Q_0}}(\cH_{Q,W}))
\]
where we write $\overline{\Module}^{\zeta}_{Q,W}(\ff,-)=\bigoplus_{\dd\in\bN^{Q_0}}\Module^{\zeta}_{Q,W}(\ff,\dd)\otimes\bL^{-\ff\cdot\dd/2}$.  Via the commutation relation $T_{\infty}\cdot T^{\dd}=q^{\ff\cdot \dd}T^{\dd}\cdot T_{\infty}$ and the ansatz \eqref{wt_BPS_def} defining the refined BPS invariants, we arrive at the equation
\[
\chi_{\wt,\bN^{Q_0}}(\overline{\Module}^{\zeta}_{Q,W}(\ff,-))=\Exp\left(\sum_{\dd\in\bN^{Q_0}\setminus \{0\}} -q^{1/2}\Omega^{\wt}_{Q,W,\dd}(q^{1/2})(1-q^{-\ff\cdot\dd})/(1-q)T^{\dd}\right).
\]
Since each of the spaces $\cM^{\ff\framed}_{\dd}(Q)$ are schemes, the vanishing cycle cohomology of each one is concentrated in bounded degrees, and so we can take Euler characteristics, i.e. we can evaluate the generating series $\chi_{\wt,\bN^{Q_0}}(\overline{\Module}^{\zeta}_{Q,W}(\ff,-))$ at $q^{1/2}=1$.    We arrive at the following infinite Euler product expansion
\[
\chi_{\bN^{Q_0}}(\overline{\Module}^{\zeta}_{Q,W}(\ff,-))=\prod_{\dd\in\bN^{Q_0}\setminus\{0\}}(1-T^{\dd})^{ \omega_{Q,W,\dd} \ff\cdot \dd}
\]
where $\omega_{Q,W,\dd}\coloneqq \Omega_{Q,W,\dd}(1)$ are the \textit{unrefined} BPS invariants.  The above derivation leads to the following \textit{definition} of unrefined BPS invariants, or Donaldson--Thomas invariants, which avoids the use of stacks: write $\chi_{\bN^{Q_0}}(\overline{\Module}^{\zeta}_{Q,W}(\ff,-))=\prod_{\dd\in\bN^{Q_0}}(1-T^{\dd})^{a_\dd}$ for integers $a_{\dd}$ and define $\omega_{Q,W,\dd}=a_{\dd}/(\ff\cdot \dd)$.  It is of course unclear that the unrefined BPS invariants defined this way are integers.  This integrality is implied by the cohomological integrality theorem, as shown above.

The general principal is the same in more general settings (for example when defining numerical DT invariants of categories of semistable coherent sheaves on Calabi--Yau threefolds as in \cite{MR1818182,MR2951762}): we first define a partition function $\cZ(T)$ by taking modified Euler characteristics of certain moduli \textit{schemes}, and express $\cZ(T)$ as an infinite product, and define DT invariants to be appropriate fractions of the exponents, where ``appropriateness'' (e.g. division by $\ff\cdot\dd$ in the above example) is determined by a wall-crossing argument with relation to a stacky moduli problem.  The resulting fractions are \textit{integers} if the refined BPS invariants for the stacky moduli problem are Laurent polynomials.

\section{BPS cohomology and geometric representation theory}
\label{GRT_section}
Since the right-hand side of the isomorphism $\Phi$ in Theorem \ref{int_thm} is the underlying space of an algebra, it is natural to wonder if $\cH_A$ carries an interesting algebra structure.  Indeed it does, and it is often considerably more interesting than a free supercommutative algebra; $\Phi$ is more akin to a Poincar\'e--Birkhoff--Witt (PBW) type isomorphism\footnote{This isomorphism tells us that, at the level of vector spaces, the morphism $\Sym(\fg)\rightarrow \UEA(\fg)$ to the universal enveloping algebra of a Lie algebra $\fg$, given by evaluating symmetric tensors, is an isomorphism.} than an actual isomorphism of algebras.  To define the \textit{BPS Lie algebra} that plays a role similar to that of the Lie algebra $\fg$ in the PBW theorem, we delve into the theory of cohomological Hall algebras (CoHAs).
\subsection{Cohomological Hall induction}
\label{CoHI}
Fix a quiver $Q$ with potential $W$, and define $A=\Jac(Q,W)$.  As above, we allow the case $W=0$, in which case there is an isomorphism $A\cong \bC Q$.  Given dimension vectors $\dd',\dd''\in \bN^{Q_0}$ we set $\dd=\dd'+\dd''$ and define $\bA_{\dd',\dd''}(Q)\subset \bA_{\dd}(Q)$ to be the subspace corresponding to $\dd$-dimensional $\bC Q$-modules preserving the $Q_0$-graded partial flag $0\subset \bigoplus_{i\in Q_0} \bC^{\dd'_i}\subset \bigoplus_{i\in Q_0}\bC^{\dd_i}$.  This is acted on by the product of parabolic groups $\GL_{\dd',\dd''}\subset \GL_{\dd}$ preserving the same flag.  We set $\fM_{\dd',\dd''}(Q)=\bA_{\dd',\dd''}(Q)/\GL_{\dd',\dd''}$, which is equivalent to the stack of short exact sequences $0\rightarrow \rho'\rightarrow \rho\rightarrow \rho''\rightarrow 0$ of $\bC Q$-modules of dimension vector $\dd',\dd,\dd''$ respectively.  There are morphisms $\pi_1,\pi_2,\pi_3$ from $\fM_{\dd',\dd''}(Q)$ to the stacks $\fM_{\dd'}(Q)$, $\fM_{\dd}(Q)$ and $\fM_{\dd''}(Q)$ respectively, which at the level of points takes such a short exact sequence to $\rho',\rho$ or $\rho''$ respectively.

The pullback morphism $(\pi_1\times \pi_3)^{\star}\colon \HO(\fM_{\dd'}(Q),\bQ)\otimes \HO(\fM_{\dd''}(Q),\bQ)\rightarrow \HO(\fM_{\dd',\dd''}(Q),\bQ)$ is defined using the K\"unneth isomorphism and functoriality of cohomology.  The pullback morphism $\HO(\fM_{\dd'}(Q),\phi_{\Tr(W)}\bQ)\otimes \HO(\fM_{\dd''}(Q),\phi_{\Tr(W)}\bQ)\rightarrow \HO(\fM_{\dd',\dd''}(Q),\phi_{\Tr(W)}\bQ)$ is defined the same way: now the Thom--Sebastiani isomorphism \cite{MR1818986} plays the role of the K\"unneth isomorphism.  The pushforward morphism $\HO(\fM_{\dd',\dd''}(Q),\bQ)\rightarrow \HO(\fM_{\dd}(Q),\bQ)$ is defined via the Verdier dual $\pi_{2,*}\bQ_{\fM_{\dd',\dd''}(Q)}\otimes\bL ^{-\dim}\rightarrow \bQ_{\fM_{\dd}(Q)}\otimes\bL^{-\dim'}$ of the adjunction morphism, $\bQ_{\fM_{\dd}(Q)}\rightarrow \pi_{2,*}\bQ_{\fM_{\dd',\dd''}(Q)}$, using the projectivity of the morphism $\pi_2$. There is an overall jump in cohomological degree/Tate twist, given by the difference $\dim-\dim'$ between the dimensions of the source and target of $\pi_2$.  Using commutativity of vanishing cycles functors along proper morphisms we obtain the morphism $\pi_{2,*}\phi_{\Tr(W)}\bQ_{\fM_{\dd',\dd''}(Q)}\otimes\bL^{-\dim}\rightarrow \phi_{\Tr(W)}\bQ_{\fM_{\dd}(Q)}\otimes\bL^{-\dim'}$, and the pushforward morphism in vanishing cycle cohomology by taking derived global sections.  For details of this construction we refer to \cite{MR2851153,MR3667216}.

Combining the above two operations we obtain a morphism
\[
\pi_{\star}\circ (\pi_1\times\pi_3)^{\star}\colon \cH_{A,\dd'}\otimes \cH_{A,\dd''}\rightarrow \cH_{A,\dd}
\]
which \textit{does} preserve the cohomological degree if $Q$ is symmetric, and is moreover a morphism of mixed Hodge structures.  Taking the direct sum across all $\dd',\dd''\in\bN^{Q_0}$ we define the operation $\mult\colon \cH_A\otimes\cH_A\rightarrow \cH_A$.  It is proved in \cite[Section 7]{MR2851153} that this endows $\cH_{A}$ with the structure of an associative algebra: the \textit{Kontsevich--Soibelman CoHA}.

\begin{remark}
An extra structure on this algebra comes from the \textit{determinant line bundle}.  For each $\dd\in\bN^{Q_0}$ the stack $\fM_{\dd}(A)$ carries a canonical rank $\sum_{i\in Q_0}\dd_i$ vector bundle $V_{\mathrm{taut}}$ given by the underling vector space of the tautological family of $\dd$-dimensional $A$-modules, and we set $\Det=\wedge^{\mathrm{top}}V_{\mathrm{taut}}$.  This yields a morphism $\fM_{\dd}(A)\rightarrow \pt/\bC^*$, and, taking cohomology, an action of $\bQ[u]=\HO(\pt/\bC^*,\bQ)$ on $\cH_A$ by derivations.  Here $u$ acts by multiplication by $c_1(\Det)$.  In the language of \cite{BDINKP}, $\Det$ provides a \textit{global equivariant parameter}.
\end{remark}

\begin{example}
Set $W=0$.  In this case, the mixed Hodge structure on $\cH_A$ is pure, of Tate type, and we ignore it, considering $\cH_Q$ as just a $\bZ_{\coh}\oplus\bN^{Q_0}$-graded vector space.  Then (after a sign-twist of the usual monoidal structure on the category of $\bZ_{\coh}\oplus\bN^{Q_0}$-graded vector spaces) it is shown in \cite{MR2851153} that the algebra $\cH_Q$ is a supercommutative algebra, and the isomorphisms of Examples \ref{qdlog_example} and \ref{boson_example} are isomorphisms of algebras\footnote{The same is true of Examples \ref{olqp_example} and \ref{olqp_example}: proving this is an instructive exercise.  Hint: show that in these examples, the perverse filtration introduced below agrees with the cohomological filtration.}.
\end{example}

The construction above is in fact a very special case of the \textit{cohomological Hall induction} constructed by Kinjo, Park and Safronov in \cite{KPS}.  In [ibid], the authors demonstrate how to give the mixed Hodge structure $\HO(\fM^{\sst}_{\sA},\phi_{\DT})$ an associative algebra structure for each of the examples (1), (2), (3) in \S \ref{general_BPS}, by solving a much more general version of this problem: the cohomological Hall induction problem for DT sheaves, for correspondence diagrams given by graded/filtered points in the framework developed by Halpern--Leistner.  See \cite{KPS, Hal14} for details.

\subsection{PBW theorem and BPS Lie algebra}
\label{BPS_Lie_section}
Fix a symmetric quiver $Q$ with potential $W$ and set $A=\Jac(Q,W)$.  We define the perverse filtration on $\cH_A^{\crit}$ as in \S \ref{perv_filt_section}.  The action of $c_1(\Det)$ sends $\fP_j\cH_A$ to $\fP_{j+2}\cH_A$.  The $\bQ[u]$-action restricts to a morphism $\fn^+_A\otimes\bQ[u]\hookrightarrow \cH_A$ which is injective.
\begin{theorem}[\cite{MR4132957}]
\label{PBW_thoerem}
The multiplication $\mult$ respects the perverse filtration on $\cH_A$.  The associated graded algebra $\Gr_{\fP}(\cH_A)$ is a free supercommutative algebra, and the morphism $\Sym\left(\fn_A^+\otimes\bQ[u]\right)\rightarrow \cH_A$ built from the algebra structure on $\cH_{A}$ is an isomorphism of $\bZ_{\coh}\oplus\bN^{Q_0}$-graded mixed Hodge structures.
\end{theorem}
By the first part of the theorem, the commutator Lie bracket $[-,-]$ is zero in the associated graded algebra $\Gr_{\fP}(\cH_A)$, so the commutator bracket applied to $\fP_1\cH_A\otimes\fP_1\cH_A$ has image in $\fP_1\cH_A=\fn_A^+$.  In other words, the mixed Hodge structure $\fn^+_A\subset \cH_A$ is closed under the commutator Lie bracket, and since this is precisely the subspace encoding the refined BPS invariants for the category of $A$-modules, we call this the \textit{BPS Lie algebra}.

There is a relative version of Theorem \ref{PBW_thoerem}.  Recall that the direct sum morphism $+\colon \cM_A\times\cM_A\rightarrow \cM_A$ makes the category of (locally bounded below) complexes of mixed Hodge modules (or perverse sheaves) on $\cM_A$ into a symmetric tensor category.  Then $\cA_A\coloneqq \bigoplus_{\dd\in\bN^{Q_0}}\JH_{\dd,*}\phi_{\Tr(W)}\bQ\otimes\bL^{\chi_Q(\dd,\dd)/2}$ is an algebra object in this category, meaning that there is a morphism $\cA_A\boxdot\cA_A\rightarrow \cA_A$ satisfying the usual axioms of an associative algebra.  Then \cite[Theorem C]{MR4132957} says that the isomorphism of Theorem \ref{rel_int_thm} can be realised by evaluating inside this algebra object.  By construction, this relative PBW theorem is a special case of the statement that \eqref{genint} is an isomorphism.
\subsection{CoHAs of sheaves with zero-dimensional support}
\label{zds_sheaves}
Let $Q=Q^{(1)}$ be the Jordan quiver from Example \ref{boson_example}.  The tripled quiver $\tilde{Q}$ defined in \S \ref{Kac_section} has three loops $x,x^*,\omega$, and $\tilde{W}=[x,x^*]\omega$.  Setting $A=\Jac(\tilde{Q},\tilde{W})$ we have $A\cong \bC[x,x^*,\omega]$.  By dimensional reduction there is an isomorphism of $\bZ_{\coh}\oplus\bN$-graded mixed Hodge structures $\cH^{\crit}_{\bC[x,x^*,\omega]}\cong \cH^{\BoMo}_{\bC[x,x^*]}$.  The object $\cH_{\bC[x,x^*]}$ was given the structure of a $\bZ_{\coh}\oplus\bN$-graded associative algebra by Schiffmann and Vasserot in \cite{MR3150250}, and one can show \cite{MR3727563,MR4091073} that in fact (after some mild twisting of signs) the dimensional reduction isomorphism takes the Kontsevich--Soibelman product on $\cH^{\crit}_{\bC[x,x^*,\omega]}$ to the Schiffmann--Vasserot product on $\cH^{\BoMo}_{\bC[x,x^*]}$.

Let $\mW^+_{1+\infty}=\Span_{\bQ}(z^mD^a\;\lvert\; m\in\bZ_{\geq 1},a\in\bZ_{\geq 0})$ be the span (over $\bQ$) of the indicated differential operators on $\bC^*$, where $D=z(d/dz)$.  This is a Lie algebra with Lie bracket $[z^mD^a,z^nD^b]=z^{m+n}((D+n)^aD^b-D^a(D+m)^b)$.  We give $\mW_{1+\infty}^+$ an ascending $\bZ$-filtration by setting $F_iW^+_{1+\infty}=\Span_{\bQ}(z^mD^a\;\lvert\; m\in\bZ_{\geq 1},0\leq a\leq (i+2)/2)$.  Set $\hat{\fg}=\Gr_F(W^+_{1+\infty})$.  Then there is an isomorphism of algebras $\cH^{\crit}_{\bC[x,x^*,\omega]}\cong \UEA(\hat{\fg})$ and so an isomorphism of algebras $\cH^{\BoMo}_{\bC[x,x^*]}\cong \UEA(\hat{\fg})$ via dimensional reduction.  See \cite{A2CoHA} and references therein, or \cite{MMSV} for the generalisation below.

Via the construction of \cite{MR4662292}, for \textit{any} smooth quasiprojective surface $S$ one can associate a $\bN$-graded CoHA $\cH^{\BoMo}_S$ for which the $n$th component is the Borel--Moore homology of the stack of sheaves with zero-dimensional support and length $n$.  Assuming that the cohomology of $S$ is itself pure\footnote{It is expected that this condition can be dropped without changing any of the statements that follow.}, this CoHA has been completely described in \cite{MMSV}.  If $S$ has trivial canonical bundle (as in the case of $\bA^2$) the description simplifies: define the Lie algebra $\fw^+(S)=\hat{\fg}\otimes \HO(S,\bQ)$ via the Lie bracket $[g\otimes x,g'\otimes x']=[g,g']\otimes xx'$.  Then there is an isomorphism of algebras $\cH^{\BoMo}_S\cong \UEA(\fw^+(S))$.  If $\omega_S$ is not trivial the description of $\cH^{\BoMo}_S$ is given in \cite{MMSV}; in this case it is a \textit{deformed} universal enveloping algebra.

In light of this description, one is led to the following conjecture of Minets regarding the critical CoHA of zero-dimensional sheaves on 3-Calabi--Yau varieties (defined via the Kinjo--Park--Safronov theorem \cite{KPS}):
\begin{conjecture} 
Let $X$ be a quasiprojective 3-Calabi--Yau variety.  Let $\mW^+(X)$ be the algebra generated by symbols $T^+_m(\lambda)$ with $m\in\bN$ and $\lambda\in\HO(X,\bQ)$, and subject to the relations $T^+_m(x\lambda)=xT^+_m(\lambda)$ for $x\in\mathbb{Q}$  and
\begin{align*}
&[T^+_m(\lambda\mu),T^+_n(\nu)]=[T^+_m(\lambda),T^+_n(\nu\mu)];\\
&[T^+_m(\lambda),T^+_{n+3}(\mu)]-3[T^+_{m+1}(\lambda),T^+_{n+2}(\mu)]+3[T^+_{m+2}(\lambda),T^+_{n+1}(\mu)]-[T^+_{m+3}(\lambda),T^+_{n}(\mu)]\\
&-[T^+_m(\lambda),T^+_{n+1}(c_2(X)\mu)]-[T^+_{m+1}(\lambda),T^+_{n}(c_2(X)\mu)]+\{T^+_m,T^+_n\}(\Delta_X \lambda\mu)=0;\\
&\sum_{\pi\in\fS_3}\pi\cdot [T^+_{m_3}(\lambda_3),[T^+_{m_2}(\lambda_2),T^+_{m_1+1}(\lambda_1)]]=0
\end{align*}
where $\fS_3$ is the symmetric group, and $\{T^+_m,T^+_n\}(\Delta_X \lambda\mu)$ is as in \cite[Sec.3.1]{MMSV}.  Then  $\cH_X^{\crit}\cong \mW^+(X)$.
\end{conjecture}
We may consider instead the \textit{equivariant} CoHA $\cH_{\bA^3}^{\crit,T}$, where $T\subset (\bC^*)^3$ is the rank two torus preserving the Calabi--Yau 3-form $dx\wedge dy\wedge dz$.  The conjecture predicts that $\cH_{\bA^3}^{\crit,T}\cong \Yang^+_{t_1,t_2,t_3}(\widehat{\fgl(1)})^+$, where the right-hand side is the positive half of the affine Yangian.  This is indeed true, see e.g. \cite{A2CoHA} and references therein.
\subsection{Modules from stable framed representations}
\label{mod_constr}
We introduce a class of representations of CoHAs built from \textit{stable framed representations}.  See e.g. \cite{MR3618057} for a fuller account.

Let $Q$ be a quiver, and let $Q^+$ be obtained from $Q$ by adding a new vertex $\infty$ and a finite number of arrows going from $\infty$ to $Q_0$, and from $Q_0$ to $\infty$.  The pair $(Q_{\ff},Q)$ from \S \ref{why_integrality} is a special case of this setup, in which \textit{all} of the new arrows go from $\infty$ to $Q_0$.  We fix the stability condition $\zeta\in\bQ^{Q^+_0}$ by setting $\zeta_i=0$ for $i\in Q_0$ and $\zeta_{\infty}=1$.  We define $(\dd,n)\in\bN^{Q^+_0}$ as in \S \ref{why_integrality}.  A $(\dd,1)$-dimensional $\bC Q^+$-module $\rho$ is $\zeta$-semistable if and only if $\bC Q^+ \cdot (1_{\infty}\cdot \rho)=\rho$, i.e. $\rho$ is generated by the one-dimensional vector space $1_{\infty}\cdot \rho$ as a $\bC Q^+$-module.  If $\rho\rightarrow \rho''$ is a surjection of $\bC Q^+$-modules with $\dim_{Q^+_0}(\rho)=(\dd,1)$, $\dim_{Q^+_0}(\rho'')=(\dd'',1)$ for some $\dd,\dd''\in\bN^{Q_0}$, and $\rho$ is $\zeta$-semistable, then $\rho''$ is as well.  Let $\cM^{\zeta\frex}_{\dd',\dd''}(Q)$ be the stack of pairs of such a surjection, along with a trivialisation\footnote{In the description of this stack as a global quotient stack, including the data of this trivialization has the same effect as \text{not} including the component $\bC^*$ of the gauge group $\GL_{(\dd,1)}$ acting by change of basis at the vertex $\infty$.  Compare with the construction of moduli \textit{schemes} in \S \ref{why_integrality}.} $1_{\infty}\cdot \rho\cong \bC$, where we have set $\dd'=\dd-\dd''$.  Equivalently, we define $\cM^{\zeta\frex}_{\dd',\dd''}(Q)$ to be the stack of short exact sequences of $\bC Q^+$-modules
\[
0\rightarrow \rho'\rightarrow \rho\rightarrow \rho''\rightarrow 0
\]
of dimension vectors $(\dd',0), (\dd,1), (\dd'',1)$ respectively, where $\rho$ and $\rho''$ are $\zeta$-semistable and $1_{\infty}\cdot\rho\cong\bC$ is trivialised.  We define $\cM^{\zeta\frex}_{\dd''}(Q)$ to be the moduli stack of pairs of a $\zeta$-semistable $(\dd'',1)$-dimensional $\bC Q^+$-module $\rho''$, along with a trivialisation $1_{\infty}\cdot \rho''$.  This stack can be written as a global quotient $\bA_{(\dd,1)}^{\zeta\sst}(Q^+)/\GL_{\dd}$ as in \S \ref{why_integrality}, and is in fact a smooth moduli \textit{scheme}.  We consider the correspondence diagram
\[
\fM_{\dd'}(Q)\times \cM^{\zeta\frex}_{\dd''}(Q)\xleftarrow{\pi_1\times \pi_3}\cM^{\zeta\frex}_{\dd',\dd''}(Q)\xrightarrow{\pi_2} \cM^{\zeta\frex}_{\dd}(Q)
\]
analogous to the diagram used to define the Kontsevich--Soibelman Hall product.  The crucial property of the diagram that we used in the definition of the pushforward of vanishing cycle cohomology was that $\pi_2$ was proper.  The same is true here due to the above-mentioned fact about surjections $\rho\rightarrow \rho''$. 

Let $W^+$ be a potential for $Q^+$, i.e. a linear combination of cyclic paths in $Q^+$.  We obtain a potential $W$ for $Q$ by throwing away all cyclic paths that are not entirely contained in $Q$.  We define 
\[
\Module=\bigoplus_{\dd\in\bN^{Q_0}}\HO(\cM^{\zeta\frex}_{\dd}(Q),\phi_{\Tr(W^+)}\bQ)\otimes \bL^{\chi_Q(\dd,\dd)-\dd\cdot\ff}
\]
where $\ff_i$ is the number of arrows from $\infty$ to $i\in Q_0$.  Mimicking the construction of the Kontsevich--Soibelman CoHA, the composition $\pi_{2,\star}\circ(\pi_1\times\pi_3)^{\star}$ yields a morphism $\cH_{Q,W}\otimes \Module\rightarrow \Module$, making $\Module$ into a module over $\cH_{Q,W}$.

\subsection{Nakajima quiver varieties}
\label{NQV_sec}
Let $Q$ be a quiver, let $\ff$ be a dimension vector, and consider the quiver $\widetilde{Q_{\ff}}$ obtained by framing and then tripling $Q$, as in \S \ref{why_integrality} and \S \ref{Kac_section}, and let $\widetilde{W_{\ff}}$ be the canonical cubic potential for this quiver.  Now remove the loop $\omega_{\infty}$, and denote by $W^+$ the potential obtained by removing all cyclic paths from $\widetilde{W_{\ff}}$ that pass through $\omega_{\infty}$.  We consider $\Module_{\ff, \dd}(Q)=\HO(\cM^{\zeta\frex}_{\dd}(\tilde{Q}),\phi_{\Tr(W^+)}\bQ)\otimes \bL^{\chi_Q(\dd,\dd)-\dd\cdot\ff}$ as above, noting that $\tilde{Q}$ is the quiver we get from $Q^+$ by removing $\infty$ and all arrows incident to it.

Define the preprojective algebra $\Pi_{Q_{\ff}}$ as in \S \ref{Kac_section}, and denote by $\cN_{Q}(\ff,\dd)\coloneqq \cM^{\zeta\sst}_{(\dd,1)}(\Pi_{Q_{\ff}})$ the fine moduli space of $\zeta$-semistable $(\dd,1)$-dimensional modules.  This moduli space is smooth and symplectic, and is an example of a \textit{Nakajima quiver variety} \cite{MR1302318}.  The following isomorphism from \cite[Section 6.2]{MR4661532} is obtained by dimensional reduction:
\[
\Module_{\ff, \dd}(Q)\cong \HO(\cM^{\zeta\sst}_{(\dd,1)}(\Pi_{Q_{\ff}}),\bQ)\otimes \bL^{-\dim/2}
\]
where $\dim$ is the dimension of the smooth variety $\cM^{\zeta\sst}_{(\dd,1)}(\Pi_{Q_{\ff}})$.  Note that the restriction of $W^+$ to $\tilde{Q}$ is the canonical cubic potential defined in \S \ref{Kac_section}.  So, as a special case of the construction from \S \ref{mod_constr}, we find that the cohomology of Nakajima quiver varieties forms modules over the the CoHA $\cH_{\tilde{Q},\tilde{W}}$, i.e. for every $\ff\in\bN^{Q_0}$ the direct sum $\Module_{\ff}(Q)\coloneqq \bigoplus_{\dd\in\bN^{Q_0}}\Module_{\ff, \dd}(Q)$ forms a module for the CoHA $\cH_{\tilde{Q},\tilde{W}}$.

\section{BPS cohomology for 2-Calabi--Yau categories}
\label{2CY_section}
Let $\sC$ be a 2-Calabi--Yau category.  We assume also that we have chosen an abelian heart $\sA$ of this dg category.  We will also restrict our attention to an Abelian subcategory $\sA^{\sst}$ of semistable objects, in order to be able to define a \textit{good moduli space} $p\colon \fM^{\sst}_{\sA}\rightarrow \cM_{\sA}$ in the sense of \cite{MR3237451}.  For the exact scope of Theorem \ref{DHSM2_thm} below we refer to \cite{DHSM2}, for now we just list a few examples.
\begin{enumerate}
\item
Let $Q$ be a quiver\footnote{Some care has to be taken if $Q$ is of Dynkin type.  In this case, one picks $\sC$ to be the dg category of modules over the \textit{derived} preprojective alkebra.  But in \textit{any} case the abelian heart $\sA$ is the usual category of modules over the usual preprojective algebra.}, and define $\sA=\Pi_Q\lmod$ to be the category of finite-dimensional modules over $\Pi_Q$.  We can restrict to a subcategory of semistable modules of fixed slope or consider all modules.
\item
Let $C$ be a smooth projective curve, and let $\sA^{\sst}=\Higgs_{\theta}^{\sst}(C)$ be the category of semistable Higgs sheaves of fixed slope $\theta\in\bQ\cup\{\infty\}$  on $C$.
\item
Let $\Sigma_g$ be a genus $g$ Riemann surface, and let $\sA$ be the category of finite-dimensional $\bC[\pi_1(\Sigma_g)]$-modules.  We choose the vacuous stability condition: all modules are semistable, stable modules are the simple ones.
\item
Let $S$ be a K3 or Abelian surface.  Let $\sA=\Coh(S)$.  Fix a polarization $H$ of $S$, and a monic polynomial $p(t)\in\bQ[t]$ of degree at most 2.  The objects of $\sA^{\sst}$ are the Gieseker-semistable objects with normalised Hilbert polynomial equal to $p(t)$.
\end{enumerate}

We define $N^+=\KK^{\eff}(\sA)$, the semigroup of effective classes in the numerical Grothendieck group of $\sA$, and set $N$ to be the numerical Grothendieck group of $\sA$, equipped with the symmetric bilinear form $\chi(-,-)$ given by the Euler pairing $\chi(x,y)=\sum_{i\in \bZ}(-1)^i\dim(\Ext_{\sC}^i(x,y))$.  In the example $\sA=\Pi_Q\lmod$, we have $N^+=\bN^{Q_0}$, $N=\bZ^{Q_0}$ and $\chi(x,y)=(\dim_{Q_0}(x),\dim_{Q_0}(y))_Q$.  If $\alpha\in N^+$ we denote the corresponding union of connected components of $\fM^{\sst}_{\sA}$ and $\cM_{\sA}$ by $\fM^{\sst}_{\sA,\alpha}$ and $\cM_{\sA,\alpha}$, respectively.  

We define $\cH_{\sA}\coloneqq \bigoplus_{\alpha\in N^+} \HO^{\BoMo}(\fM^{\sst}_{\sA,\alpha},\bQ)\otimes \bL^{\vdim /2}$.  There is a uniform construction making $\cH_{\sA}$ into an associative $\bZ_{\coh}\oplus N^+$-graded algebra; see \cite{MR3150250,MR4087863,MR3521588,MR4662292}.  For the category $\sA=\Pi_Q\lmod$ this is the dimensional reduction of the algebra considered in \S \ref{GRT_section}.

We define $\Sigma_{\sA}^+$, the set of \textit{primitive positive roots}, to be the set of $\alpha$ for which the subvariety $\cM^{\simp}_{\sA,\alpha}\subset \cM_{\sA,\alpha}$ of simple objects is non-empty.  Note that in the category of semistable coherent sheaves of fixed normalised Hilbert polynomial, or semistable Higgs bundles of fixed slope, simple objects are precisely the stable ones.  These roots are called real, isotropic, or hyperbolic depending on whether $\dim(\cM^{\simp}_{\sA,\alpha})=0,2,d>2$ respectively.  We define the \textit{simple positive roots} $\Phi_{\sA}^+$ to be the union of $\Sigma_{\sA}^+$ with multiples of all the isotropic roots.  In the case $\sA=\Pi_Q\lmod$ we have $N=\bN^{Q_0}$, the Euler form $\chi(-,-)$ is identified with $(-,-)_Q$, and the definitions of all of the various kinds of roots become purely combinatorial \cite[Section 5.1]{DHSM2}.  One reason that the more abstract definition makes sense for more general 2CY categories is that they all locally ``look like'' the category of $\Pi_Q$-modules: see \cite{MR4794593}.

\subsection{Generalised Kac--Moody (GKM) Lie algebras}
We pick one of the 2CY categories above and continue to set $N^+=\KK^{\eff}(\sA)$.  Let $V$ be a $\bZ_{\coh}\oplus N^+$-graded vector space.  We assume that $V_{\alpha}=\bQ$ (concentrated in cohomological degree zero) if $\alpha$ is a real simple root, and $V_{\alpha}=0$ for $\alpha\notin\Phi_{\sA}^+$.  Set $\fh=N\otimes_{\bZ}\bQ$, which we give $\bZ_{\coh}\oplus N$-degree $(0,0)$.  
We define $\fg_V$ to be the free Lie algebra generated by $V\oplus \fh\oplus V^{\vee}$ modulo the relations
\begin{itemize}
\item
$[h,h']=0$ for $h,h'\in \fh$.
\item
$[h,g]=\chi(h,\alpha) g$ for $\alpha\in N$, if $g\in (V\oplus V^{\vee})_{\alpha}$ and $h\in\fh$.
\item
$[v,v']=v'(v)\alpha\in\fh$ for $v\in V_{\alpha}$ and $v'\in V^{\vee}$.
\item
(The Serre relations): $[g,-]^{1-(\alpha,\alpha')_Q}(g')=0$ if either $(\alpha,\alpha')_Q=0$ or $\alpha$ is a real simple root.
\end{itemize}
There is a triangular decomposition $\fg_V=\fn^-_V\oplus \fh\oplus \fn^+_V$ with $\fn_V^{\pm}$ in $N$-degrees $\pm N^+\setminus \{0\}$.  Note that $\fn_V^-$ and $\fn^+_V$ are the free Lie algebras generated by $V^{\vee}$ (respectively $V$) subject only to the Serre relations.  GKM Lie algebras appear naturally in the theory of crystal graphs associated to quivers with loops; see \cite{MR2553376, MR3569998,MR4708147}.

The convolution tensor product on $\Perv(\cM_{\sA})$ makes it into a $N^+$-graded tensor category, meaning that every object $\cF$ decomposes as a direct sum $\cF=\bigoplus_{\alpha\in N^+}\cF_{\alpha}$, and $(\cF\boxdot\cG)_{\alpha}=\bigoplus_{\alpha'+\alpha''=\alpha}\cF_{\alpha'}\boxdot\cG_{\alpha''}$.  Although there are no \textit{negatively} graded objects, with respect to the $N$-grading, we can still define the analogue of the Lie algebra object $\fn^+_V$ in this category.  In place of $V$, we start with a perverse sheaf $\cG\in \Perv(\cM_{\sA})$ satisfying the property that $\cG_{\alpha}=\bQ_{\cM_{\sA,\alpha}}$ if $\alpha$ is a simple root (note that this implies that $\cM_{\sA,\alpha}\cong\pt$), and $\cG_{\alpha}=0$ if $\alpha\notin\Phi^+_{\sA}$.  We define $\fn^+_{\cG}$ to be the free Lie algebra object in the symmetric tensor category $\Perv(\cM_{\sA})$ subject to the Serre relations above; see \cite[Section 3.5]{DHSM2} for details of the construction.

If we set $V=\HO(\cM_{\sA},\cG)$ then there is an isomorphism $\HO(\cM_{\sA},\fn^+_{\cG})\cong \fn^+_V$.  In particular, the derived global sections of $\fn^+_{\cG}$ carry a Lie algebra structure, and there is a natural inclusion of Lie algebras $\HO(\cM_{\sA},\fn^+_{\cG})\hookrightarrow \fg_V$ into a doubled Lie algebra, which is a GKM Lie algebra.

\subsection{BPS Lie algebras for 2CY categories}
\begin{theorem}[\cite{DHSM2}]
\label{DHSM2_thm}
Let $\sC,\sA,\sA^{\sst}$ be as at the start of \S \ref{2CY_section}.  We define $\cG$ as follows 
\[
\cG_{\alpha}=\begin{cases} \IC_{\cM_{\sA,\alpha}}& \textrm{if there exists a simple object in }\sA^{\sst}\textrm{ of class }\alpha\\
+_{n,*}\IC_{\cM_{\sA,\alpha'}}&\textrm{if }\alpha=n\alpha' \textrm{ with }\alpha'\in\Sigma_{\sA}^+\textrm{ isotropic}\\
0&\textrm{otherwise}
\end{cases}
\]
where $+_{n}\colon \cM_{\sA,\alpha'}\rightarrow \cM_{\sA,n\alpha'}$ is the morphism taking $x$ to $x^{\oplus n}$.  Then there is a PBW isomorphism
\[
\Phi\colon \Sym_{\boxdot}\left(\fn^+_{\cG}\otimes\HO(\pt/\bC^*,\bQ)\right)\rightarrow \bigoplus_{\alpha\in N^+}p_*\bD\bQ_{\fM_{\sA,\alpha}}\otimes\bL^{\vdim/2}
\]
and taking global sections of $\Phi$ we deduce that $\cH_{\sA}$ satisfies cohomological integrality: there is a PBW isomorphism $\Sym\left(\fn^+_V\otimes\HO(\pt/\bC^*,\bQ)\right)\rightarrow \cH_{\sA}$ with the BPS cohomology given by $\fn^+_V\otimes\bL^{-1/2}$.
\end{theorem}
It follows that the BPS Lie algebra $\fn^+_V$ is the positive half of a GKM Lie algebra, which we denote $\fg_{\sA}$.
\label{DHSM_results}
\subsection{Maulik--Okounkov Yangians}
We restrict our attention to $\Pi_Q\lmod$.  As in \S \ref{zds_sheaves} we consider an equivariant theory: we consider the $T=(\bC^*)^2$-action on $\Pi_Q$ for which the first $\bC^*$-factor scales all of the arrows $a\in Q_1$ and the second factor scales all of the arrows $a^*$.  We write $\HO_T\coloneqq \HO(\pt/T,\bQ)$.  As in \S \ref{zds_sheaves} we consider the equivariant CoHA $\cH_{\Pi_Q}^T$.  Define $\Module^T_{\ff,\dd}$ to be the $T$-equivariant cohomology of the Nakajima quiver variety $\cN_Q(\ff,\dd)$.  As in \S \ref{NQV_sec} the direct sums $\Module_{\ff}^T(Q)\cong \bigoplus_{\dd\in\bN^{Q_0}}\Module^T_{\ff,\dd}(Q)$ form modules over the equivariant CoHA $\cH^T_{\Pi_Q}$.  One may identify the raising operators on cohomology of Nakajima quiver varieties considered by Nakajima \cite{MR1604167}, Grojnowski \cite{MR1386846} and Varagnolo \cite{MR1818101} with the action of explicit elements of $\cH^T_{\Pi_Q}$, although for general $Q$ the algebra $\cH^T_{\Pi_Q}$ contains many more raising operators besides these ones.

In \cite{MR3951025} Maulik and Okounkov defined a new Yangian-type algebra $\Yang^{\MO}_Q$, acting by raising (and lowering) operators on each $\Module^T_{\ff}(Q)$, defined via the theory of stable envelopes and the resulting R matrices.  As in the case of the Yangians associated to simple Lie algebras, one has a PBW-type isomorphism $\Phi:\Sym_{\HO_T}(\fg^{\MO}_Q\otimes\bQ[u])\cong \Yang^{\MO}_Q$.  Here $\fg^{\MO}_Q$ is a (free) $\HO_T$-module, so that $\fg^{\MO}\otimes\bQ[u]$ carries a $\HO_T$-module structure, and the domain of $\Phi$ is the free algebra object (in the tensor category of $\HO_T$-modules with tensor structure $\otimes_{\HO_T}$) generated by this module.  As for (generalised) Kac--Moody Lie algebras, there is a triangular decomposition $\fg^{\MO}_Q\cong \fn^{\MO,-}_Q\oplus\fh\oplus \fn^{\MO,+}_Q$.  We define $\Yang^{\MO,+}_Q$ to be the image of $\Sym_{\HO_T}(\fn^{\MO,+}_Q\otimes\bQ[u])$ under $\Phi$.  It is a subalgebra of $\Yang^{\MO}_Q$.

In \cite{BD23} we proved that there is an isomorphism $\psi\colon \fg^{T}_{\Pi_Q\lmod}\cong \fg^{T,\MO}_Q$.  It then follows from \cite[Corollary 11.9]{DHSM2} that $\fg^{T,\MO}_Q$ can actually be defined over $\bQ$: there is a Lie algebra $\fg^{\MO}_Q$ defined over $\bQ$ and an isomorphism of Lie algebras  $\fg^{T,\MO}_Q\cong \fg^{\MO}_Q\otimes \HO_T$, where $[g\otimes \lambda,g'\otimes \lambda']=[g,g']\otimes \lambda\lambda'$.  Furthermore, from \S \ref{DHSM_results} we obtain a complete description of $\fg^{\MO}_Q$ as a $\bZ_{\coh}\oplus\bZ^{Q_0}$-graded GKM algebra.  From Theorem \ref{Kac_to_BPS} we deduce the following equality, conjectured by Okounkov, and proven independently in \cite{SVOC}: $\poinc(\fg_{Q,\dd}^{\MO})=\kac_{Q,\dd}(q^{-1})$.  Moreover, $\Psi$ extends to an isomorphism $\cH_{\Pi_Q}^T\cong \Yang^{\MO,+}_Q$ \cite{BD23,SVOC}, so that we can fully describe the Yangian via a full description of $\cH_{\Pi_Q}^T$.  This strategy is completed for the case of cyclic quivers by Jindal \cite{ShivangYangian}.
\subsection{Nonabelian Hodge isomorphisms}
\label{NAHT_section}

Fix a smooth projective curve $C$, and let $\Sigma_g$ be the underlying Riemann surface.  To start with we consider semistable Higgs bundles of slope $\theta=0$; in order to work with other slopes, on the Betti side of the nonabelian Hodge correspondence one should consider \textit{twisted} representations of $\pi_1(\Sigma_g)$ below, i.e. representations of the punctured Riemann surface $\Sigma_g\setminus\{p\}$ for which the monodromy around $p$ is multiplication by a prescribed root of unity.  We write $\fM^{\Dol}_{r,d}$ for the moduli stack of rank $r$ and degree $d$ semistable Higgs bundles, $\cM^{\Dol}_{r,d}$ for the coarse moduli space, and $\fM^{\Betti}_r$ and $\cM^{\Betti}_r$ for the moduli stack and coarse moduli space of $r$-dimensional $\pi_1(\Sigma_g)$-representations.

We write $\cH^{\Betti}=\HO^{\BoMo}(\fM^{\Betti}_{\bullet},\bQ)\otimes \bL^{\vdim/2}$ and $\cH^{\Dol}=\HO^{\BoMo}(\fM^{\Dol}_{\bullet,0},\bQ)\otimes \bL^{\vdim/2}$.  By Theorem \ref{DHSM2_thm} these two $\bZ_{\coh}\oplus\bN$-graded vector spaces are determined by the intersection cohomology of $\cM^{\Betti}$ and $\cM^{\Dol}$ respectively.  These spaces are homeomorphic, via the classical nonabelian Hodge correspondence, and it follows that there is an isomorphism \cite{DHSM1} of graded vector spaces $\Lambda\colon\cH^{\Betti}\cong \cH^{\Dol}$ (see \cite{Henn23} for an upgrade to an algebra isomorphism).

This isomorphism does \textit{not} lift to an isomorphism of mixed Hodge structures.  The $\chi$-independence theorem of \cite{KiKo21} states that, considered as filtered vector spaces, with the filtration defined via perverse filtration functors with respect to the direct image to the base of the Hitchin system, there is an isomorphism $\BPS^{\Dol}_{r}\cong \HO(\cM_{r,1}^{\Dol},\bQ)[-2(1-g)r^2]$.  It follows that there is an isomorphism of cohomologically graded vector spaces $\BPS^{\Betti}_{r}\cong \HO(\cM^{\Betti,\twisted}_r,\bQ)[-2(1-g)r^2]$ where the target is the cohomology of the twisted character variety.  This cohomology provides one side of the P=W conjecture relating this cohomology with that of $\cM_{r,1}^{\Dol}$.  The P=W conjecture is now a theorem \cite{HMMS, MR4792069,MR4860586}.  Putting all this together, it is natural to make the following:
\begin{conjecture}
$\Lambda$ takes the weight filtration on $\cH^{\Betti}$ to a (doubled) perverse filtration on $\cH^{\Dol}$.
\end{conjecture}
We refer to \cite{MR4543450} for the precise definition of the perverse filtration on $\cH^{\Dol}$, the precise form of the conjecture, and the relation to the PI=WI conjecture of de Cataldo and Maulik \cite{MR4221002}.

\section{BPS cohomology beyond moduli stacks of objects}

Many stacks of interest are not moduli stacks of objects in a linear category.  For instance, let $M$ be a real oriented closed three-manifold without boundary, and consider the stack $\fM_G(M)$ of $G$-local systems for $G$ a connected reductive group not equal to $\mathrm{GL}_n$.  By \cite{PTVV}, this stack has a (-1)-shifted symplectic structure, and by \cite{NS23} it carries orientation data, so that we can define the DT sheaf $\phi_{\DT}$ as in \S \ref{general_BPS}.  It is not immediately clear how to \textit{define} the BPS sheaf or integrality for $\fM_G(M)$.   For example if $G$ is semisimple, the decomposition of $p_*\phi_{\DT}$ into cohomologically shifted mixed Hodge modules starts in degree \textit{zero}, not one, where $p\colon \fM_G(M)\rightarrow \cM_G(M)$ is the good moduli space morphism.

Nonetheless, there \textit{is} a version of BPS cohomology and integrality for this stack, which is special case of a general integrality theorem for (-1)-shifted stacks; we refer the reader to \cite{BDINKP,HeKi25} for all the details, and also the analogous results for $(0)$-shifted and smooth stacks.  For a finite type (oriented) (-1)-shifted symplectic stack $\fM$ with good moduli space $p\colon \fM\rightarrow \cM$, the analogue of Proposition \ref{vc_decomp_thm} is:
\begin{proposition}[\cite{BDINKP} Proposition 7.2.9]
There is an isomorphism $p_*\phi_{\DT}\cong \bigoplus_{i\geq \mathrm{crk}(\fM)} \cH^i\left(p_*\phi_{\DT}\right)[-i]$ in the derived category of unbounded complexes of mixed Hodge modules on $\cM$.
\end{proposition}
Here $\crk(\fM)$ is the maximal $r$ such that there exists an action of $\pt/T$ on $\fM$ with $T$ a rank $r$ torus, that does not factor through an action of $\pt/T'$ for $T'$ of lower rank.  The proposition motivates the definitions 
\[
\shBPS_{\fM}\coloneqq\cH^{0}(p_*\phi_{\DT}\otimes\bL^{-\mathrm{crk}(\fM)/2});\quad\quad\BPS_{\fM}\coloneqq \HO(\cM,\shBPS_{\fM}).
\]
Write $c_G=\mathrm{crk}(\fM_G(M))$.  The cohomological integrality isomorphism for $\fM_G(M)$ then states that there is an isomorphism $\bigoplus_{L\subset G}\left( \BPS_{\fM_{L}(M)}\otimes \bL^{c_L/2}\otimes \HO(\pt/(\bC^*)^{c_L},\bQ)\right)^{W_L}\cong \HO(\fM_G(M),\phi_{\DT})$ where $W_L$ is the relative Weyl group of $L$, and the sum is over conjugacy classes of Levi subgroups of $G$.  We refer to \cite[Section 10.3]{BDINKP} for details.  The existence of this isomorphism was proved in \cite{BDINKP} under a certain orthogonality assumption on the tangent weights of $\fM_G(M)$, and then proved without this assumption in \cite{HeKi25}.  

We expect that the study of BPS cohomology for stacks of local systems on three manifolds will give a powerful toolkit for defining and understanding enumerative invariants in low-dimensional topology.  For example, we finish with the following conjecture of Pavel Safronov, linking Betti Langlands duality to the study of this generalised BPS cohomology:
\begin{conjecture}
For $(G,G^{\vee})$ Langlands dual groups, there is an isomorphism $\BPS_{\fM_{G}(M)}\cong \BPS_{\fM_{G^{\vee}}(M)}$.
\end{conjecture}
For $M=(S^1)^3$ the conjecture was proved by Kaubrys \cite{Kaub24} for $G=\mathrm{SL}_p$ with $p$ prime, and is proved for many more Langlands dual pairs of groups in forthcoming work of Hennecart and Kinjo.


\section*{Acknowledgments.}
I was supported by a Royal Society University Research Fellowship (no. 221040) during the writing of this article.

I have many people to thank for their enlightening work, conversations, and collaboration over the years regarding this subject, including Tommaso Botta, Tom Bridgeland, Jim Bryan, Chenjing Bu, Francesca Carocci, Kevin Costello, Tudor Dimofte, Camilla Felisetti, Ezra Getzler, Victor Ginzburg, Tam\'as Hausel, Lucien Hennecart, Victoria Hoskins, Andr\'es Ib\'a\~nez N\'u\~nez, Shivang Jindal, David Jordan, Dominic Joyce, \u{S}ar\=unas Kaubrys, Bernhard Keller,  Tasuki Kinjo, Maxim Kontsevich, Kentaro Nagao, Naoki Koseki, Travis Mandel, Davesh Maulik, Mirko Mauri, Michael McBreen, Anton Mellit, Sven Meinhardt, Sasha Minets, Andrei Negu\c{t}, Andrei Okounkov, Tudor P\u{a}durariu, Andrea Ricolfi, Pavel Safronov, Francesco Sala, Olivier Schiffmann,  Sebastian Schlegel Mejia, Junliang Shen, Vivek Shende, Yan Soibelman, Bal\'azs Szendr\H{o}i, Richard Thomas, Yukinobu Toda, Okke van Garderen, Michel van Garrel, Eric Vasserot, Michael Wemyss and Dimitri Wyss.

\bibliographystyle{siamplain}
\bibliography{BD_Bibliography}
\end{document}